\documentclass[12pt,seceqn,secthm,amsthm]{elsart}
\usepackage{amsfonts}
\usepackage{amsmath}
\usepackage{amssymb}
\usepackage[dvips]{graphicx}%

\theoremstyle{definition}
\newtheorem{definitions}[thm]{Definitions}
\newtheorem{notation}[thm]{Notation}
\theoremstyle{remark}
\newtheorem{remarks}[thm]{Remarks}
\newcommand{\newsection}[1]{\setcounter{equation}{0} \setcounter{theorem}{0}
\section{#1}}
\renewcommand{\newsection}{\section}
\renewcommand{\theenumi}{\alph{enumi}}

\newcounter{enumd}

\newcounter{enumdd}

\newcommand{\NI}{\noindent}
\newcommand{\IR}{I\!\!R}
\renewcommand{\IR}{\mathbb{R}}

\newcommand{\IC}{I\!\!C}
\renewcommand{\IC}{\mathbb{C}}
\newcommand{\IN}{I\!\!N}
\renewcommand{\IN}{\mathbb{N}}
\newcommand{\IT}{I\!\!T}
\renewcommand{\IT}{\mathbb{T}}
\newcommand{\clk}{{\cal K}}
\renewcommand{\clk}{\mathcal{K}}
\newcommand{\clh}{{\cal H}}
\renewcommand{\clh}{\mathcal{H}}
\newcommand{\clo}{{\cal O}}
\renewcommand{\clo}{\mathcal{O}}
\newcommand{\clm}{{\cal M}}
\renewcommand{\clm}{\mathcal{M}}
\renewcommand{\phi}{\varphi}
\newcommand{\raro}{\rightarrow}
\newcommand{\be}{\begin{equation}}
\newcommand{\ee}{\end{equation}}
\newcommand{\Imultiindex}{I}
\newcommand{\Iidentity}{\hbox{\upshape \small1\kern-3.3pt\normalsize1}}
\newcommand{\Iinterval}{\mathcal{I}}

\journal{J. Approx. Theory}

\begin{document}
\begin{frontmatter}

\title{Localized bases in $L^{2}(0,1)$ and their use in the analysis of Brownian motion}

\author{Palle E. T. Jorgensen\corauthref{cor1}\thanksref{petj}}\thanks
[petj]{This material is based upon work supported by the U.S.
National Science Foundation under grants DMS-0139473 (FRG) and DMS-0457581.}
\address{Department of Mathematics,
The University of Iowa,
14 MacLean Hall,
Iowa City, IA 52242-1419,
U.S.A.}
\ead{jorgen@math.uiowa.edu}
\ead[url]{http://www.math.uiowa.edu/\symbol{126}jorgen}
\author{Anilesh Mohari}
\address{S. N. Bose Centre for Basic Sciences,
JD Block, Sector-3, Calcutta-98,
India}
\ead{anilesh@boson.bose.res.in}
\corauth[cor1]{Tel: (319) 335-0782, Fax: (319) 335-0627, E-mail: jorgen@math.uiowa.edu}

\begin{abstract}
Motivated by
problems on Brownian motion, we 
introduce a recursive scheme for a basis construction in the Hilbert space
$L^{2}(0,1)$ which is analogous to that of Haar and Walsh.
More generally, we find a new decomposition theory
for the Hilbert space of square-integrable
functions on the unit-interval, both with respect to Lebesgue measure, and
also with respect to a wider class of self-similar measures $\mu$. That is, we
consider recursive and orthogonal decompositions for the Hilbert space
$L^2(\mu)$ where $\mu$ is some self-similar measure on $[0, 1]$.
Up to two specific
reflection symmetries, our scheme produces infinite families of orthonormal
bases in $L^{2}(0,1)$. Our approach is as versatile as the more traditional
spline constructions. But while singly generated spline bases typically do not
produce orthonormal bases, each of our present algorithms does.

\end{abstract}

\begin{keyword}
Haar, Walsh, orthonormal basis, Hilbert space, Cuntz relations,
irreducible representation, wavelets, iterated function system, Cantor set
\renewcommand{\MSC}{{\par\leavevmode\hbox{\it2000 MSC:\ }}}\MSC
42C40, 37F40, 46E22, 47L30
\end{keyword}

\end{frontmatter}

\newsection{\label{Int}Introduction}

The basis constructions considered in this paper involve elements from the 
theory of operator algebras. Since this may not be widely known to readers in 
approximation theory, we begin with a few remarks. 

The study of operator 
algebras breaks up in two parts: One the study of ``the algebras themselves" 
as they emerge from the axioms, von Neumann algebras, and $C^*$-algebras. The 
other has a more applied slant: It involves ``representations" of the 
algebras. There is a close connection between the two parts of the theory: For 
example, representations of $C^*$-algebras generate von Neumann algebras. It 
was realized in the last ten years (see e.g., \cite{BJP96}, \cite{DuJo05b}) 
that certain families of representations are useful in basis constructions in 
harmonic analysis, in approximations, in signal/image analysis, and more 
generally in computational mathematics. The bases in question may typically be 
built up from representations of an especially important family of simple 
$C^*$-algebras, known as the Cuntz algebras \cite{Cu77}. These Cuntz algebras 
(see Lemma \ref{LEMMA2.2} below) are denoted 
$\mathcal{O}_{2}, \mathcal{O}_{3},..,$ including $\mathcal{O}_{\infty}$.

The connection to Cuntz algebras $\mathcal{O}_{N}$ is further relevant to the kind of 
dynamical systems built on iterated branching-laws, with the case of 
$\mathcal{O}_{N}$ 
representing $N$-fold branching. The reason for this is that if $N$ is fixed, 
$\mathcal{O}_{N}$ includes in its definition an iterated branching, taking the 
form of 
subdivision, but now within the context of Hilbert space; so we generate 
subdivisions into orthogonal families of subspaces of the initial Hilbert 
space. In this paper, we follow up on a certain probabilistic aspect of this 
construction.

       The $\mathcal{O}_{N}$  point of view is especially well suited to basis 
constructions in such contexts as wavelets and fractals since they naturally 
involve the same kind of sub-division. Our starting point is an initial 
Hilbert space $\mathcal{H}$, where $\mathcal{H}$ may be 
$L^{2}(\mathbb{R}^{d})$, or $\mathcal{H}$ may be $L^{2}(\mu)$ for some fractal 
measure $\mu$; see Section \ref{S2} below. The more successful bases in Hilbert space 
are the orthonormal bases ONBs, and we shall consider a certain computational 
algorithm for generating them. 

       A further reason the subdivision schemes in Hilbert space are useful is 
that the more familiar Fourier wave functions are periodic, and so not 
localized. Moreover these existing Fourier tools are typically not friendly to 
algorithmic computations. By a \textit{local} (See Definition \ref{D:DEF}) construction we mean an algorithm which
begins with a finite family of functions (often one or two) having a fixed 
compact (i.e., local) support, and a procedure allowing assigned scaling and 
translation operations. As is popular in the sub-band approach to wavelets and 
wavelet packets in 1D, the scaling is typically in powers of a fixed base, 
i.e., it could consist of powers $N^{j}$ where $N$ is fixed ($N > 1$) and 
where $j$ varies over $\mathbb{Z}$, stretching and squeezing the support.

                    The main result in this paper concerns bases in the
Hilbert spaces $L^{2}(X,\mu)$ defined from measures $\mu$ arising as equilibrium
measures (also called self-similar measures) for iterated function systems
(IFS). However, our results are spelled out in more detail for
$L^2((0,1);\operatorname{Lebesgue})$, where the classical Walsh system is a special case. Our
construction uses ideas from dynamical systems and operator algebras
(specifically representations of the Cuntz algebras). The Cuntz algebra
\cite{Cu77} is used in the construction of bases for geometric structures with
self-similarity such as iterated function systems (IFS). Recall that for a
fixed finite $N$, the Cuntz algebra $\clo_N$ is generated abstractly by $N$
isometries. In representations of $\clo_N$ on a concrete Hilbert space $\clh$, the
resulting isometries $S_i$, say, have orthogonal ranges which form a partition
of unity in the particular Hilbert space $\clh$ which carries the representation,
in the sense that the identity operator $\Iidentity _{\clh}$ is written as a sum
of the $N$ projections $S_i S_i^*$ onto the respective ranges $S_i \clh$. Since the
subdivision process can be iterated, this idea has already proved useful in
understanding orthogonal families in Hilbert spaces built on IFSs; see, e.g.,
\cite{DuJo03,JoPe96,JoPe98}. Our analysis here uses such particular $\clo_N$
representations in combination with certain graph-theoretic considerations.
In addition to the reference \cite{Cu77}, the Cuntz algebras $\clo_N$ and
their representations are reviewed in \cite[Sections 7.6 and 7.7]{Jor06}. The
book \cite{Jor06} also includes additional motivation, and \cite[Chapter 4, p.~69, and
Section 9.4]{Jor06} cover details for an IFS-family of Cantor systems and their
self-similar measures.

Applications to Brownian motions are intended but postponed to a later paper.
The connection between Brownian motion, Cuntz algebras, and IFSs is treated
in the literature, for example in \cite[pp.~56--57, 151, and 203]{Jor06}.

There has been a recent increased interest in basis constructions outside the
traditional setting of harmonic analysis. The setting which so far has proved
more amenable to an explicit analysis with basis functions involves a mix of
analysis and dynamics, and it typically goes beyond the standard and more
familiar setting of orthonormal bases consisting of Fourier frequencies. The
context of frames in Hilbert space (see, e.g., \cite{ALTW04}, \cite{BoPa05},
and \cite{BPS03}) is a case in point.

As is well known, the classical setting of Fourier analysis presupposes a choice of
Fourier frequencies, or Fourier trigonometric basis functions. However, this
unduly limits our choices, and the applications: As is well known, Fourier's
basis functions are less localized, and the computational formulas typically
are not recursive. There are now alternative dynamical approaches which are
recursive, and at the same time are amenable to harmonic basis constructions;
see, e.g., \cite{Dut05}, \cite{JoPe98}, \cite{DuJo05b}, and \cite{DuJo05c}. Moreover, these
recursive models arise in applications exhibiting a suitable
\textit{scale-similarity}. Their consideration combines classical ideas from
infinite convolution with ideas from dynamics of a more recent vintage.

In this paper, we introduce a recursive scheme for a basis construction in the Hilbert space
$L^{2}(0,1)$ which is analogous to that of Haar \cite{Haa10} and Walsh
\cite{Wal23,Chr55}. While computationally efficient, these more
traditional approaches limit the choices of functions too much, often to step
functions, or at any rate to functions that have a limited number of derivatives.

We begin here with a certain dual system of axioms for reflection symmetries
for functions on the unit interval $\Iinterval =[0,1]$. We then show that up to this
reflection symmetry, we get recursive algorithms which in turn produce
infinite families of orthonormal bases localized in $L^{2}(0,1)$. Our scheme
is as versatile as the more traditional spline constructions; see, for example,
\cite{Wic94}. But while singly generated spline systems do not produce
orthonormal bases, each of our present algorithms does. Moreover our scheme is
adapted to the fixed unit interval $\Iinterval =[0,1]$ while the more traditional
wavelet-based wavepackets are designed for the construction
of orthonormal bases in $L^{2}(\IR)$; see \cite{CoWi92}. 
And if the starting functions are of compact support,
the size of the support reaches outside the unit interval $[0,1]$.

\begin{figure}[ptb]
\setlength{\unitlength}{0.775bp} \begin{picture}(426,642)(3,3)
\put(17,603){\makebox(0,0)[r]{$\scriptstyle\varphi_{0}$}}
\put(0,486){\includegraphics[bb=0 0 108 162,width=108\unitlength]{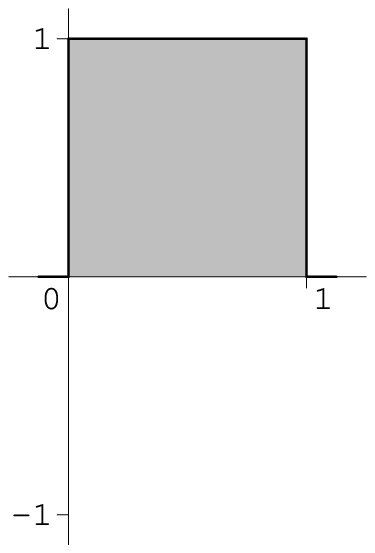}}
\put(125,603){\makebox(0,0)[r]{$\scriptstyle\varphi_{1}$}}
\put(108,486){\includegraphics[bb=0 0 108 162,width=108\unitlength]{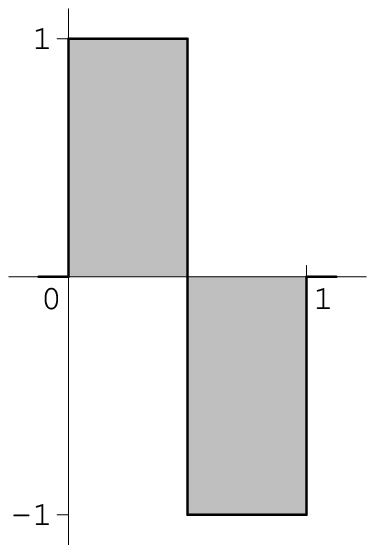}}
\put(233,603){\makebox(0,0)[r]{$\scriptstyle\varphi_{2}$}}
\put(216,486){\includegraphics[bb=0 0 108 162,width=108\unitlength]{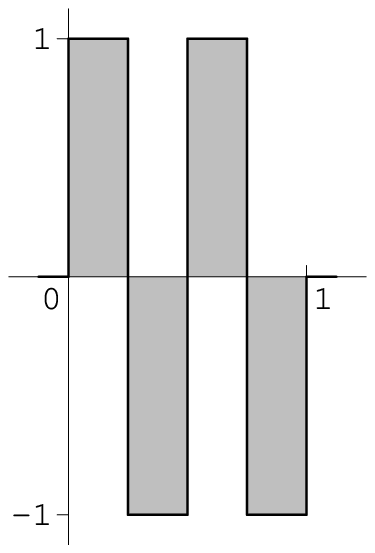}}
\put(341,603){\makebox(0,0)[r]{$\scriptstyle\varphi_{3}$}}
\put(324,486){\includegraphics[bb=0 0 108 162,width=108\unitlength]{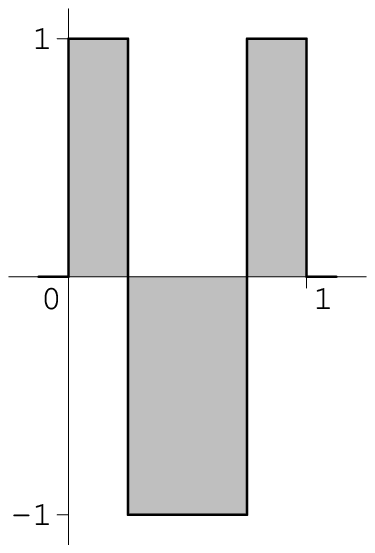}}
\put(17,441){\makebox(0,0)[r]{$\scriptstyle\varphi_{8}$}}
\put(0,324){\includegraphics[bb=0 0 108 162,width=108\unitlength]{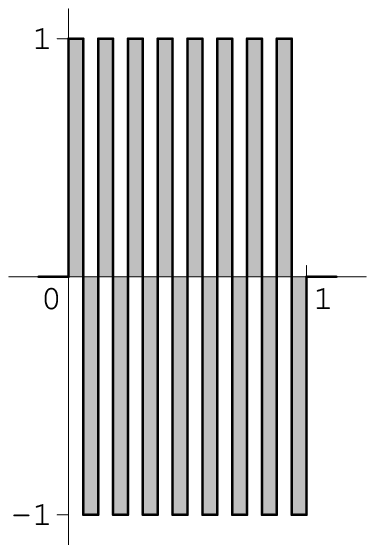}}
\put(125,441){\makebox(0,0)[r]{$\scriptstyle\varphi_{9}$}}
\put(108,324){\includegraphics[bb=0 0 108 162,width=108\unitlength]{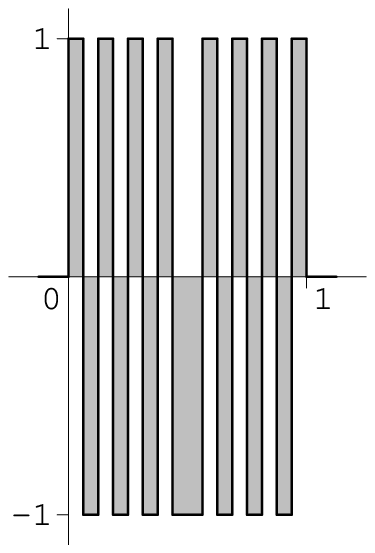}}
\put(233,441){\makebox(0,0)[r]{$\scriptstyle\varphi_{10}$}}
\put(216,324){\includegraphics[bb=0 0 108 162,width=108\unitlength]{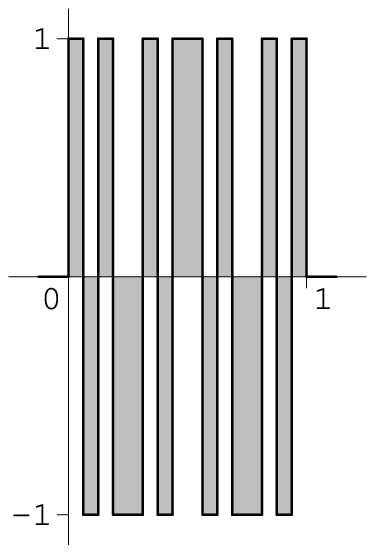}}
\put(341,441){\makebox(0,0)[r]{$\scriptstyle\varphi_{11}$}}
\put(324,324){\includegraphics[bb=0 0 108 162,width=108\unitlength]{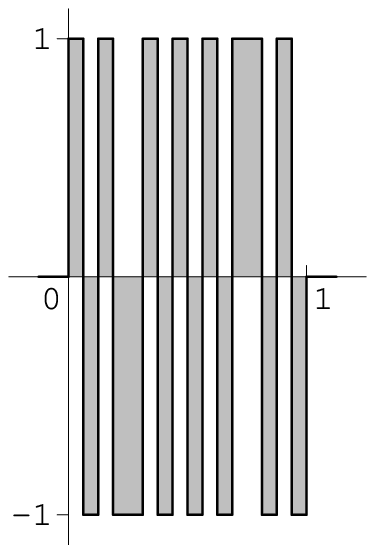}}
\put(17,279){\makebox(0,0)[r]{$\scriptstyle\varphi_{16}$}}
\put(0,162){\includegraphics[bb=0 0 108 162,width=108\unitlength]{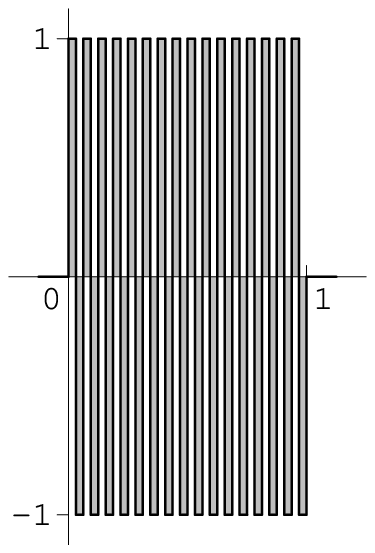}}
\put(125,279){\makebox(0,0)[r]{$\scriptstyle\varphi_{17}$}}
\put(108,162){\includegraphics[bb=0 0 108 162,width=108\unitlength]{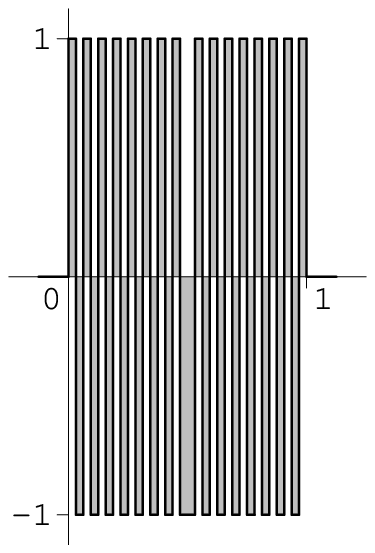}}
\put(233,279){\makebox(0,0)[r]{$\scriptstyle\varphi_{18}$}}
\put(216,162){\includegraphics[bb=0 0 108 162,width=108\unitlength]{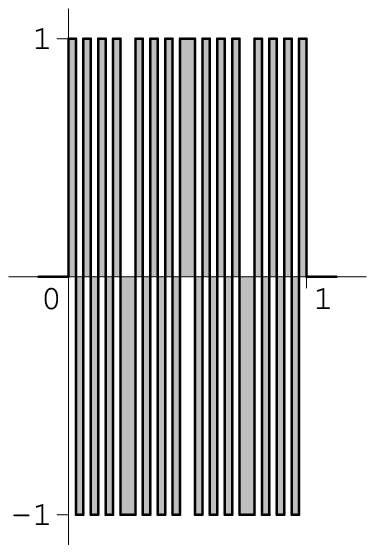}}
\put(341,279){\makebox(0,0)[r]{$\scriptstyle\varphi_{19}$}}
\put(324,162){\includegraphics[bb=0 0 108 162,width=108\unitlength]{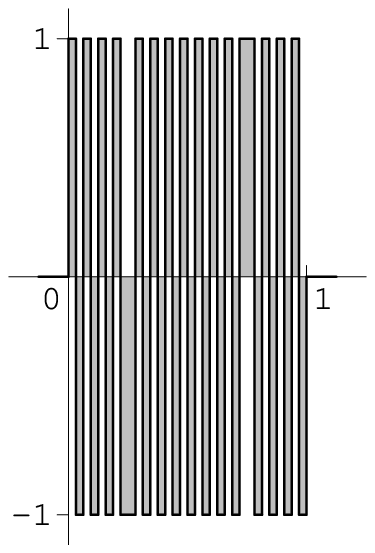}}
\put(17,117){\makebox(0,0)[r]{$\scriptstyle\varphi_{24}$}}
\put(0,0){\includegraphics[bb=0 0 108 162,width=108\unitlength]{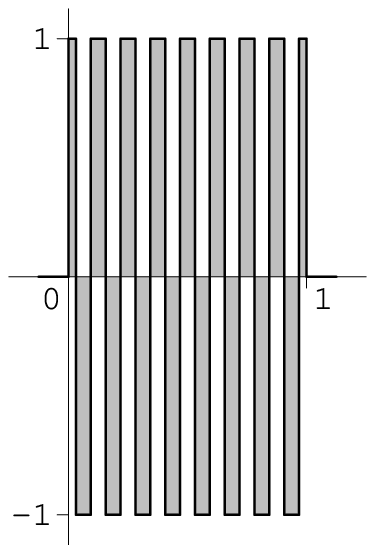}}
\put(125,117){\makebox(0,0)[r]{$\scriptstyle\varphi_{25}$}}
\put(108,0){\includegraphics[bb=0 0 108 162,width=108\unitlength]{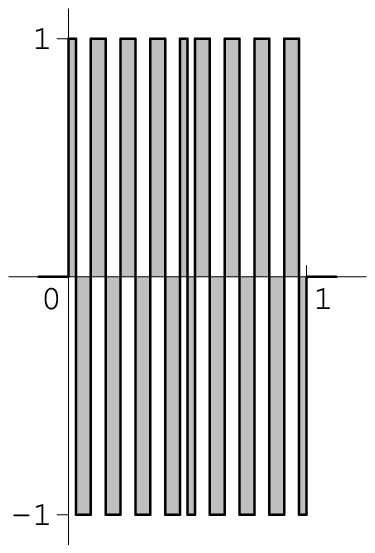}}
\put(233,117){\makebox(0,0)[r]{$\scriptstyle\varphi_{26}$}}
\put(216,0){\includegraphics[bb=0 0 108 162,width=108\unitlength]{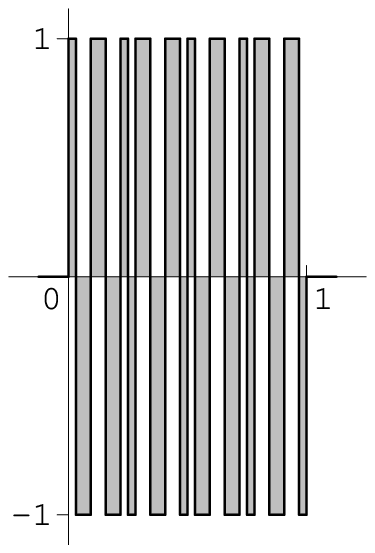}}
\put(341,117){\makebox(0,0)[r]{$\scriptstyle\varphi_{27}$}}
\put(324,0){\includegraphics[bb=0 0 108 162,width=108\unitlength]{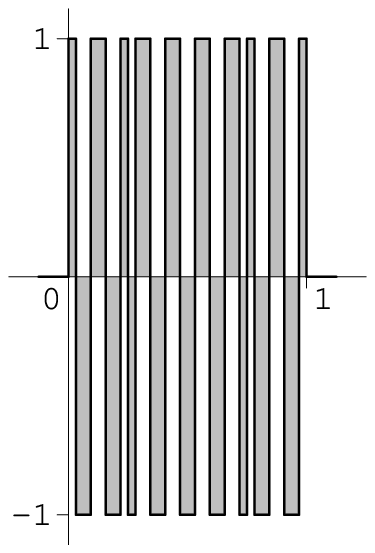}}
\end{picture}
\caption{The first thirty-two functions in the sequence $\varphi_{n}$}%
\label{FigHaarSequence}%
\end{figure}

The more traditional approaches to basis algorithms further have limited the
libraries of functions to be used at the initial step of the recursion. We get
around this here by identifying a set of symmetry conditions that may be
imposed on two functions $\psi_{0}$ and $\psi_{1}$ in the Hilbert space
$L^{2}(0,1)$;
see Fig.\ \ref{FigHaarSequence} for an illustration in the simplest
case.
Our algorithm is then based on a certain matrix scaling and
subdivision applied to these two functions. Hence our starting point is
different from the more traditional one which begins with a scaling identity,
masking coefficients, and a so-called father function $\phi_{0}$ which solves
the corresponding scaling identity; see \cite{Dau92}.

\begin{figure}[ptb]
\setlength{\unitlength}{0.775bp} \begin{picture}(426,673)(3,-28)
\put(17,603){\makebox(0,0)[r]{$\scriptstyle\varphi_{4}$}}
\put(0,486){\includegraphics[bb=0 0 108 162,width=108\unitlength]{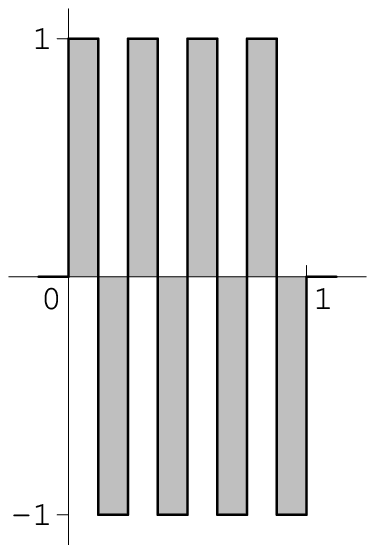}}
\put(125,603){\makebox(0,0)[r]{$\scriptstyle\varphi_{5}$}}
\put(108,486){\includegraphics[bb=0 0 108 162,width=108\unitlength]{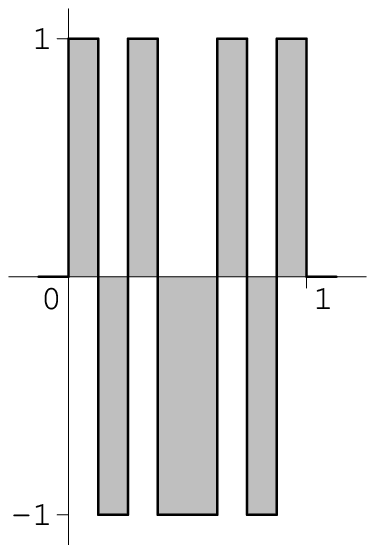}}
\put(233,603){\makebox(0,0)[r]{$\scriptstyle\varphi_{6}$}}
\put(216,486){\includegraphics[bb=0 0 108 162,width=108\unitlength]{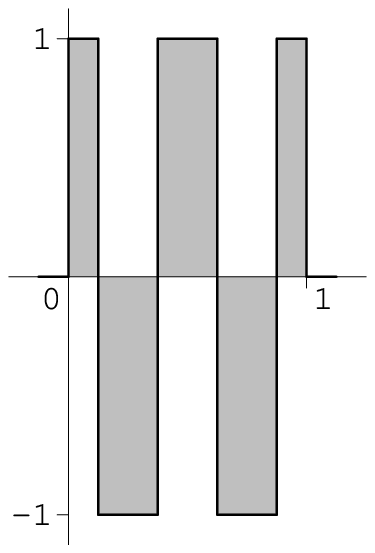}}
\put(341,603){\makebox(0,0)[r]{$\scriptstyle\varphi_{7}$}}
\put(324,486){\includegraphics[bb=0 0 108 162,width=108\unitlength]{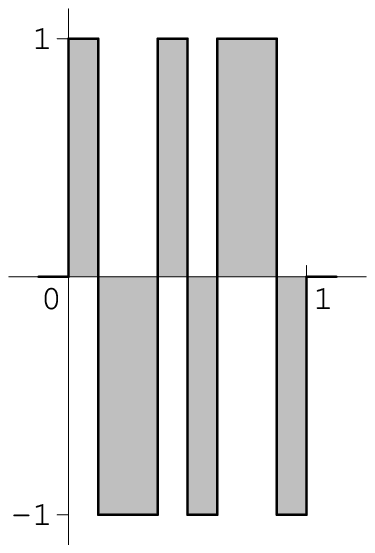}}
\put(17,441){\makebox(0,0)[r]{$\scriptstyle\varphi_{12}$}}
\put(0,324){\includegraphics[bb=0 0 108 162,width=108\unitlength]{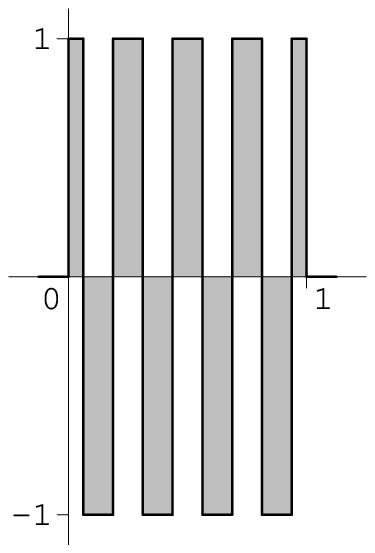}}
\put(125,441){\makebox(0,0)[r]{$\scriptstyle\varphi_{13}$}}
\put(108,324){\includegraphics[bb=0 0 108 162,width=108\unitlength]{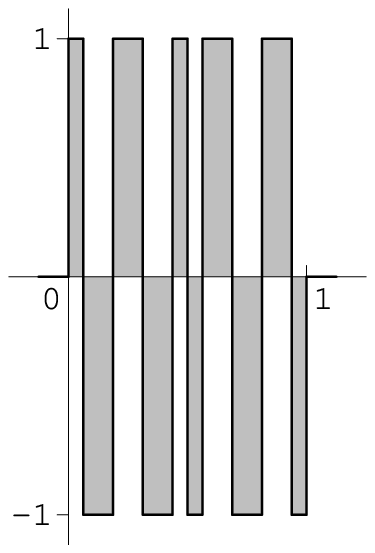}}
\put(233,441){\makebox(0,0)[r]{$\scriptstyle\varphi_{14}$}}
\put(216,324){\includegraphics[bb=0 0 108 162,width=108\unitlength]{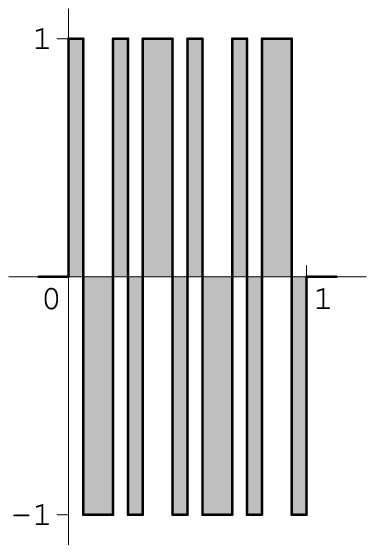}}
\put(341,441){\makebox(0,0)[r]{$\scriptstyle\varphi_{15}$}}
\put(324,324){\includegraphics[bb=0 0 108 162,width=108\unitlength]{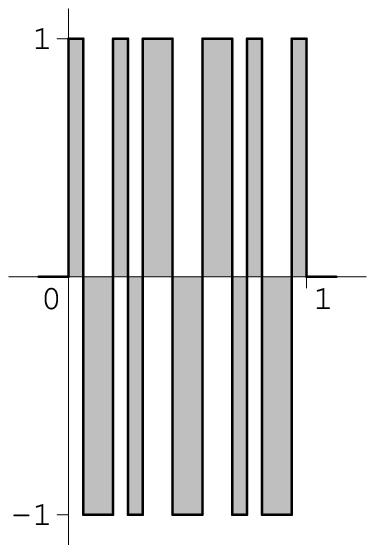}}
\put(17,279){\makebox(0,0)[r]{$\scriptstyle\varphi_{20}$}}
\put(0,162){\includegraphics[bb=0 0 108 162,width=108\unitlength]{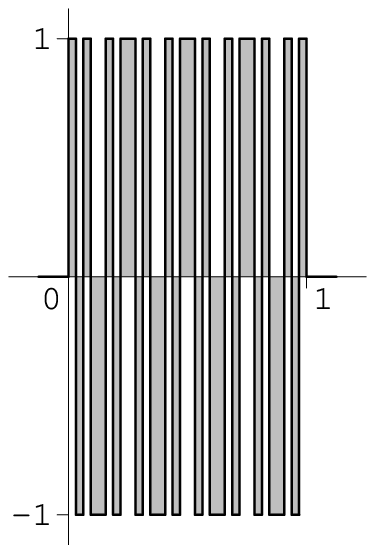}}
\put(125,279){\makebox(0,0)[r]{$\scriptstyle\varphi_{21}$}}
\put(108,162){\includegraphics[bb=0 0 108 162,width=108\unitlength]{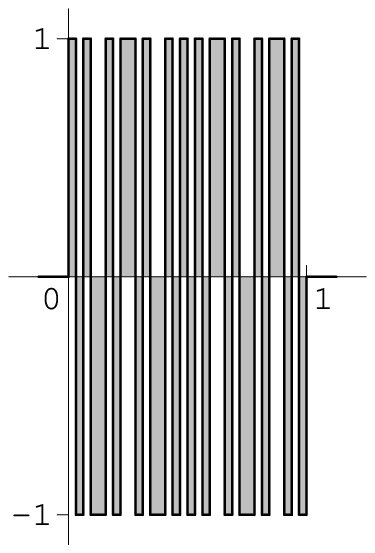}}
\put(233,279){\makebox(0,0)[r]{$\scriptstyle\varphi_{22}$}}
\put(216,162){\includegraphics[bb=0 0 108 162,width=108\unitlength]{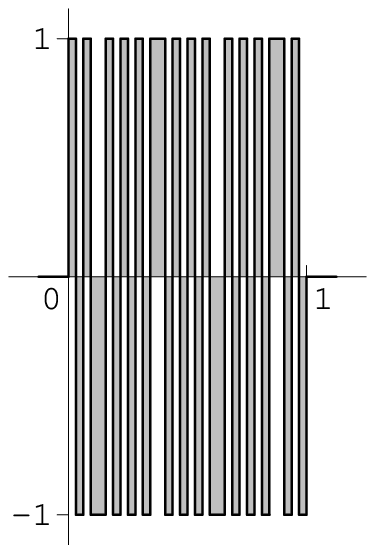}}
\put(341,279){\makebox(0,0)[r]{$\scriptstyle\varphi_{23}$}}
\put(324,162){\includegraphics[bb=0 0 108 162,width=108\unitlength]{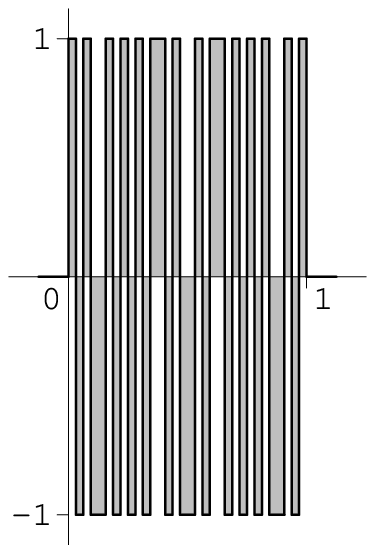}}
\put(17,117){\makebox(0,0)[r]{$\scriptstyle\varphi_{28}$}}
\put(0,0){\includegraphics[bb=0 0 108 162,width=108\unitlength]{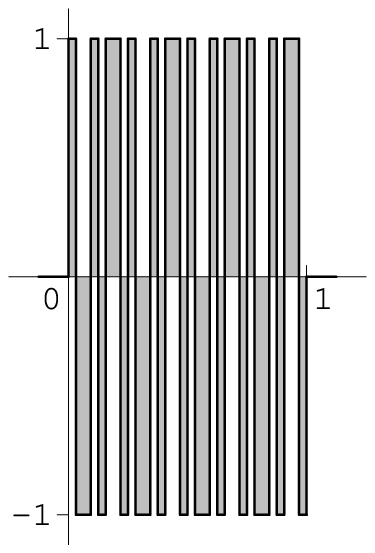}}
\put(125,117){\makebox(0,0)[r]{$\scriptstyle\varphi_{29}$}}
\put(108,0){\includegraphics[bb=0 0 108 162,width=108\unitlength]{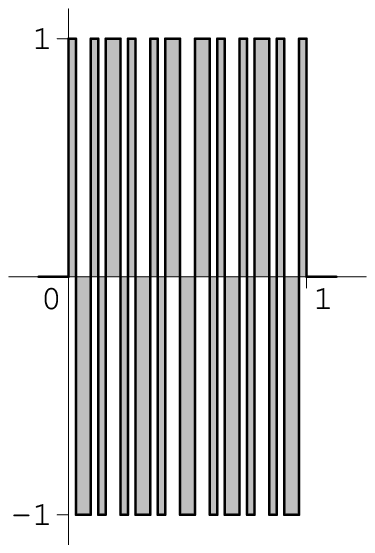}}
\put(233,117){\makebox(0,0)[r]{$\scriptstyle\varphi_{30}$}}
\put(216,0){\includegraphics[bb=0 0 108 162,width=108\unitlength]{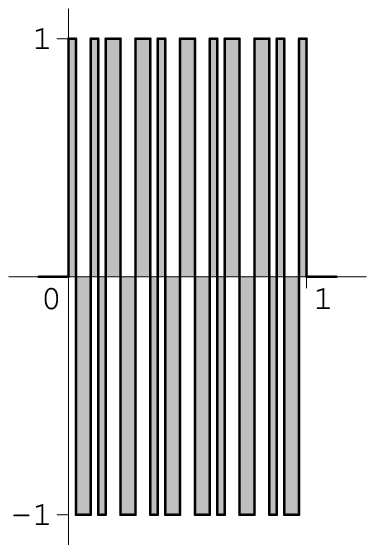}}
\put(341,117){\makebox(0,0)[r]{$\scriptstyle\varphi_{31}$}}
\put(324,0){\includegraphics[bb=0 0 108 162,width=108\unitlength]{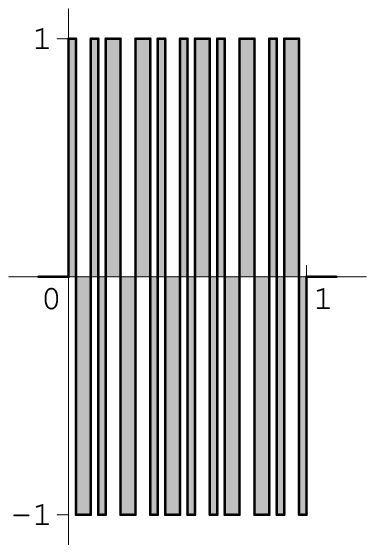}}
\put(216,-28){\makebox(0,0)[b]{\small Fig.\ \ref{FigHaarSequence} (cont'd).}}
\end{picture}%
\end{figure}%

Our justification for the term ``wavelet'' in connection with the present basis is threefold:

\NI (a) Our functions are localized in a sense which will be made clear.

\NI (b) Our construction is recursive.

\NI (c) Our algorithm for constructing orthonormal bases starts with a
prescribed and carefully selected finite system $S$ of functions. Two
operations are applied recursively to $S$, scaling and reflection. But note
that as the algorithm runs, the reflections are scaled as well. 
\begin{definitions}
\label{D:DEF}
A family of functions (typically an orthonormal basis) $\{ \psi_n: n \ge 1 \}$ 
in $L^2(0,1)$ is said to be local if for all $\epsilon > 0$ there exists an 
$N(\epsilon) \ge 1$ so that the closed linear space spanned by
$\{\psi_n: 1 \le n \le N(\epsilon) \}$ admits a finite subfamily of orthogonal 
functions $\{ \phi_n: 1 \le n \le M(\epsilon) \}$  total for the subspace 
spanned by $\{\psi_n: 1 \le n \le N(\epsilon) \}$ with 
the restriction that the Lebesgue measure of the support of each  $\phi_n$ is 
less then $\epsilon$. 
\end{definitions}

         To understand the two discrete operations which underlie our
construction, it is helpful to review a fundamental feature of wavelets in
the non-standard setting of iterated function systems (IFS); see also
\cite{JoPe94,JoPe96} and Section \ref{SNew4} below for additional details. In the simplest
of settings, i.e., that of the unit interval, the two discrete operations
going into our algorithms are that of iterated scaling by $2$, i.e., the
system $x \raro 2^m x$, modulo $1$,  for $m = 0, 1, 2, \dots$, and mid-point reflections.
Both operations have analogues for more general IFSs, and we discuss these
generalizations below.
\pagebreak

         Our use of the terms ``reflection'' and ``reflection symmetry" is
related to, but different from, the one studied in, for example, \cite{JoOl98}. The
main difference is that in our present context, the Hilbert-space inner
product is preserved after the reflection, while in \cite{JoOl98} it is changed,
i.e., it is subjected to a certain renormalization.

\newsection{\label{S2}Iterated function systems}

To understand our results it is useful to consider a slightly more general 
setup: Let d be a fixed dimension, and consider a given finite set S of affine 
and contractive mappings $\tau_{i} : \mathbb{R}^{d} \to \mathbb{R}^{d}$. 
There are several 
interesting limits in the literature, arising from iteration of such a system 
$S = (\tau_{i})$. The accepted terminology is ``affine Iterated Function 
System (IFS)". See \cite{Hut81}, \cite{Jor06}, and formulas 
(\ref{eqUniqueBorel})-(\ref{eq2.1}) below.

           Following Cantor's middle-third construction, note that one 
limiting object derived from $S$ results from recursive iterations of the 
individual maps in $S$; it is an attractor which takes the form of a compact 
subset $X (= X(S))$ of $\mathbb{R}^{d}$.  The set $X(S)$ often has fractal like 
properties: And in one dimension ($d = 1$) the deleted-middle-third Cantor set 
is an example, but the unit interval is one too.  The other iteration limits 
in this context take place in the family of probability measures on 
$\mathbb{R}^{d}$.  
However when $S$ is given, the limit measure in question (called equilibrium 
measure) depends further on a chosen assignment of probability weights: It 
turns out that, given $S$, and given a fixed assignment of probabilities 
($p_{i}$) to the $\tau_{i}$s, there is a unique equilibrium measure 
$\mu = \mu_{S,p}$. If $p_{i} > 0$ for all $i$, then the support of the 
measure $\mu_{S,p}$ is $X(S)$. Even in 1D when $X(S)$ may be the unit 
interval, there is a variety of measures arising from the second limit 
construction other than the restricted Lebesgue measure.

          In this paper, we are interested in localized orthonormal bases 
(ONBs) in $L^{2}(\mu_{S,p})$. It turns out that our main issues may be best 
presented in the case of $d = 1$, and in the special case where the weights 
are uniform, i.e., $p_{i} = 1/N$ where $N$ is the number of maps $\tau_{i}$ 
from the initial system $S$. In fact, as noted in \cite{DuJo03}, the Hilbert 
space $L^{2}(\mu_{S,p})$ does not have ONBs of complex exponential Fourier 
bases in the case of non-uniform weights.

While the gist of our paper is for the unit interval, we wish to add that our
construction in fact works in a more general context, that of iterated
function systems (IFS) \cite{Hut81,JoPe94}. However, it is easier to get an overview of the
totality of admissible bases in the special case of the unit interval
$\Iinterval =[0,1]$. 

We consider decomposition theory for the Hilbert space of square-integrable
functions on the unit interval, both with respect to Lebesgue measure, and
also with respect to a wider class of self-similar measures $\mu$ \cite{JoPe96}. That is, we
consider recursive and orthogonal decompositions for the Hilbert space
$L^2(\mu)$ where $\mu$ is some self-similar measure on $[0, 1]$.

We now turn to the technical details needed in our discussion of IFSs in
Sections \ref{SNew3} and \ref{SNew4} below.
For the axiomatics of IFSs, see, e.g.,
\cite{Hut81} and \cite{Jor05}. A finite iterated function system is
determined by a finite system of contractive endomorphisms $\tau=(\tau_{i})$
in a compact metric space $X$.
We shall consider here a finite
system consisting of $N$ endomorphisms of a compact space $X$. As we will see,
one of the endomorphisms will be singled out, and it will be convenient to
index it by zero, i.e., the first map is $\tau_0$.
When such a system $\tau$ is given, it is then
known that there is a unique Borel probability measure $\mu$ on $X$ such that
\begin{equation}
\mu=\frac{1}{N}\sum_{i}\mu \circ \tau_{i}^{-1},
\label{eqUniqueBorel}
\end{equation}
or equivalently,
$$
\int_{X}f(x)\,d\mu(x)= {1 \over N} \sum_i \int f(\tau_i(x))\,d\mu(x) \mbox{\qquad for all } f \in C(X).
$$
The measure $\mu$ is called the Hutchinson measure, and it is also called the
balanced invariant measure of $\tau$. In the case when $X$ is $[0,1]$, the
unit interval, then the two maps $\tau_{0}$ and $\tau_{1}$ may be taken to be
$x \rightarrow x/2$ and $x \rightarrow (x+1)/2$, respectively, and $\mu$ is then
the standard normalized Lebesgue measure on $[0,1]$. For general IFSs
($X,\tau$), the starting point for our ONB construction in the Hilbert space
$L^{2}(X,\mu)$ is then a specified finite set of functions ($\psi_{i}$) which
satisfy a certain reflection symmetry which we proceed to describe in detail
in 
Sections \ref{SNew3} and \ref{SNew4}
below.

\bigskip
Let $(X,d)$ be a compact metric space. Let $\sigma\colon X \rightarrow  X$ be 
an endomorphism such that the number of elements in $\sigma^{-1}\{ x \}$ is equal to $N$ for all 
$x \in X$, where $N$, $2 \le N < \infty$, is fixed. Iteration of branches 
$$
\sigma^{-1}(\{x \})= \{\,y \in X \mid \sigma(y)=x \,\}
$$
then gives rise to a combinatorial tree. If 
\[
\omega=(\omega_1,\omega_2,\dots) \in \Omega = \{ 0,1,\dots,N-1 \}^{\IN},
\]
an associated path may be thought of as an infinite extension of finite
walks 
$$\tau_{\omega_n}\tau_{\omega_{n-1}}\cdots\tau_{\omega_2}\tau_{\omega_1}x$$
with starting point $x$, where $(\tau_i)$, $i=0,1,\dots,N-1$, is a system of 
Borel measurable inverses of $\sigma$, i.e.,
$$\sigma \circ \tau_i= \operatorname{id}_{X},\qquad 0 \le i \le N-1 . $$
We assume that each $(\tau_i)$ is contractive, i.e., there exists a constant $0 < c < 1$ such that 
$$d(\tau_ix,\tau_iy) \le c d(x,y),\qquad x,y \in X,\;i \in \{0,1,\dots,N-1 \}.$$
This contractive condition ensures that there exists a Borel probability measure $\mu$ on $X$
such that 
$$ \mu = {1 \over N} \sum^{N-1}_{i=0}\mu \circ \tau_i^{-1} .$$
Before we state our first simple lemma we also set 
$$\alpha_i = \sigma \;\;\;\mbox{on}\;\; \tau_i(X).$$
For each $j \in \{ 0,1,\dots,N-1 \}$ we define a linear operator $(S_j)$ by
\be
S_jf = \sum^{N-1}_{k=0} e^{i2\pi {jk \over N}}
\chi_{\tau_{k}\left(  X\right)  }
f \circ \alpha_k
\label{eq2.1}
\ee
for all $f \in \mathcal{H}= L^2(X, \mu)$.

\begin{lem}
\label{LEMMA2.1}
If the operators $S_j$  in $L^2(X, \mu)$ are given as in \textup{(\ref{eq2.1})}, then the formula
for the corresponding adjoint operators $S^*_j$ is as follows:
\be
 S^*_jf= {1 \over N} \sum^{N-1}_{k=0}e^{-i2\pi {jk \over N}} f \circ \tau_k
 \ee
\end{lem}

\begin{pf}
We leave the details to the reader. They are based on the
$\tau_i$-equipartition property for the measure $\mu$.
See (\ref{eqUniqueBorel}), and also \cite{JoPe94}. 
\qed
\end{pf}

           In the next lemma we show that the operators $S_j$ from (\ref{eq2.1})
define a representation of the Cuntz $C^*$-algebra. As noted in \cite{JoPe94,JoPe96},
this $C^*$-algebra has found a variety of uses in iterated function
systems (IFS), and in approximation theory. Its use in analysis and in
physics was initiated in the paper \cite{Cu77}. For each IFS with $N$
endomorphisms, there is an associated representation of the Cuntz algebra
$\clo_N$ with $N$ generators $S_j$. These generators are also called the fundamental
isometries. (When we say ``isometry'', we are here referring to the Hilbert
space $L^2(X, \mu)$.)  From \cite{Cu77} we further know that specifying a
representation of $\clo_N$ is equivalent to specifying a system of $N$ fundamental
isometries.

\begin{lem}
\label{LEMMA2.2} 
The system of operators $(S_j)$ from \textup{(\ref{eq2.1})} defines a representation of the
Cuntz algebra $\clo_N$. We write
$(S_j) \in \operatorname{Rep}(\clo_d,\mathcal{H})$, i.e., 
\be
S_j^*S_k=\delta_{j,k} \Iidentity _{\mathcal{H}},\qquad\sum^{N-1}_{k=0} S_kS^*_k=\Iidentity _{\mathcal{H}} .
\label{eqCuntzRelations}
\ee
\end{lem}

\begin{pf}
Again, we leave the details to the reader. The argument uses the previous
lemma combined with the axioms for the system $(\tau_j)$ and the associated
measure $\mu$, outlined above.

Returning to formula (\ref{eq2.1}) and setting%
\[
m_{j}\left(  x\right)  =\sum_{k=0}^{N-1}e^{i2\pi\frac{jk}{N}}\chi_{\tau
_{k}\left(  X\right)  }\left(  x\right)  ,
\]
we see that (\ref{eq2.1}) may be rewritten in the form%
\begin{equation}
\left(  S_{j}f\right)  \left(  x\right)  =m_{j}\left(  x\right)  f\left(
\sigma\left(  x\right)  \right)  \text{\qquad for }f\in L^{2}\left(
X,\mu\right)  ,\;x\in X.\label{eq2.1rewritten}%
\end{equation}

Hence using an idea from \cite{BJMP05}, we then note that the Cuntz relations
(\ref{eqCuntzRelations}) for the operators $\left(  S_{j}\right)  $ in
(\ref{eq2.1rewritten}) are equivalent to the assertion that the $N\times N$
matrix function $U$ given by%
\[
U\left(  x\right)  =
\left(  U_{j,k}\left(  x\right)  \right)
_{j,k=0}^{N-1}  :=
\frac{1}{\sqrt{N}}\left(  m_{j}\left(  \tau
_{k}\left(  x\right)  \right)  \vphantom{\frac{1}{\sqrt{N}}}\!\right)
_{j,k=0}^{N-1}%
\]
takes values in the group $\mathrm{U}_{N}\left(  \mathbb{C}\right)  $ of all
$N\times N$ unitary matrices.

We now check that $UU^{\ast}=\Iidentity _{N}$. Substituting $m_{j}\left(  \tau
_{k}\left(  x\right)  \right)  =e^{i2\pi\frac{jk}{N}}$ into the formula for
the matrix product, we get%
\[
\sum_{k=0}^{N-1}U_{j,k}\left(  x\right)  \,\overline{U_{l,k}\left(  x\right)
}=\frac{1}{N}\sum_{k=0}^{N-1}e^{i2\pi\frac{\left(  j-l\right)  k}{N}}%
=\delta_{j,l}\,.
\qed
\]
\end{pf}

\begin{rem}
\label{Remark.p}There is a variety of representations of the Cuntz algebra
$\clo_N$, other than those that come naturally from IFSs (as in (\ref{eq2.1}) and the
lemma). In fact by a theorem of Glimm \cite{Gli60}, the set of equivalence
classes of irreducible representations of $\clo_N$ cannot be parametrized by a
Borel cross section.
\end{rem}

      Our focus in this paper is the search for orthonormal bases (ONBs) in
the Hilbert space $L^2(\mu)$ associated naturally with a fixed contractive IFS.
Once an ONB is chosen, it may be used in the analogue-to-digital (A-to-D)
conversion of signals. Conversely, for a given signal (read, $f \in L^2(\mu)$)
there is a variety of ONBs, and a choice must be made; for example, we may
wish to minimize the entropy (information loss), 
\[
\varepsilon^{\left(  \psi\right)  }\left(  f\right)  =-\sum_{k}\left\vert
\left\langle \,\psi_{k}\mid f\,\right\rangle \right\vert ^{2}\ln\left\vert
\left\langle \,\psi_{k}\mid f\,\right\rangle \right\vert ^{2},
\]
where $(\psi_k)$ is the chosen ONB. (Here, we shall normalize $f$, i.e., assume
$\left\Vert f \right\Vert _{L^2(\mu)} = 1$.)

      We make use of representations in two ways: (1) as a first step toward
the determination of an ONB; and (2) for entropy computations, even when an
ONB is not known.

       The details on (1) will be presented in the next section; and we now
briefly discuss (2).

        As illustrated in the lemma below, our use of the fundamental
isometries $(S_i)$ in a particular representation of $\clo_N$ in some Hilbert space
$\clh$ immediately yields scales of mutually orthogonal subspaces of $\clh$, or
equivalently, mutually orthogonal projections which constitute partition of
the identity operator $\Iidentity $ in $\clh$. Example: the projections $P_i := S_i S_i^*$ are
orthogonal and satisfy 
\[
\sum P_i = \Iidentity . 
\]
        As stressed in for example \cite{Kri05} and \cite{JKK05}, the use of subspaces
as opposed to vectors is significant for the design of quantum algorithms
such as quantum error-correction codes. The reason is that direct operations
on individual vectors in the Hilbert space $\clh$ typically would destroy the
quantum states on which the algorithm is operating. The quantum states (or
qubits) might represent polarized photons.  (Recall that in the conventions
of quantum theory, unit-vectors in the underlying Hilbert space $\clh$ represent
quantum states!) Our use of representations of $\clo_N$ lets the algorithms act
on subspaces in $\clh$, and as a result leave the individual quantum states
intact. The intrinsic orthogonality of the subspaces is what yields quantum
channels; see \cite{JKK05}.

\begin{lem}
\label{Lemma.p}Let $(S_i) \in \operatorname{Rep} (\clo_N, \clh)$ be a representation of $\clo_N$ acting on a
Hilbert space $\clh$. Let $\clm_k$ denote the set of all multi-indices $J = (j_1 \, j_2 \, \dots
\, j_k) $, where each $j_i$ is in $\left\{0, 1,\dots ,N - 1\right\}$. 
Then for each $k$, 
\[
S_{J}:=S_{j_{1}}\cdots S_{j_{k}},\qquad J\in\mathcal{M}_{k},
\]
is a representation of $\clo_{N^k}$; in particular, $P_J := S_J S_J^*$ yields a
commuting family of mutually orthogonal projections in $\clh$, i.e., 
\[
P_{J}P_{J^{\prime}}=\delta_{J,J^{\prime}}P_{J},\qquad J\in\mathcal{M}_{k},
\]
and%
\[
\sum_{J\in\mathcal{M}_{k}}P_{J}=\Iidentity _{\mathcal{H}}.
\]
For $f$ in $\clh$, $\left\Vert f \right\Vert = 1$, the entropy number 
\[
\varepsilon_{k}\left(  f\right)  :=-\sum_{J\in\mathcal{M}_{k}}\left\Vert P_{J}f\right\Vert ^{2}%
\ln\left\Vert P_{J}f\right\Vert ^{2}%
\]
satisfies 
\[
\varepsilon_{k+1}\left(  f\right)  =\varepsilon_{1}\left(  f\right)  +\sum
_{i}\left\Vert S_{i}^{\ast}f\right\Vert ^{2}\varepsilon_{k}\left(  \frac
{P_{i}f}{\left\Vert P_{i}f\right\Vert }\right)  ,
\]
where $P_i = S_i S_i^*$   and 
\[
\varepsilon_{1}\left(  f\right)  :=-\sum_{i=0}^{N-1}\left\Vert S_{i}^{\ast
}f\right\Vert ^{2}\ln\left\Vert S_{i}^{\ast}f\right\Vert ^{2}.
\]
\end{lem} 

\begin{pf}
The steps in the proof of the lemma are straightforward, and we leave
them to the reader. \qed
\end{pf}

\begin{exmp}
\label{EXAMPLE2.3}
Let $X=[0,1]$ and let $\sigma\colon X \raro X$ be the endomorphism $\sigma(x)= 2x \bmod 1$. Here 
$\tau_0(x)={x \over 2}$ and $\tau_1(x)= {x+1 \over 2} $ are two maps for which $\sigma \circ \tau_i(x)=x,\;i=0,1$. Let
$\alpha_0(x)=2x \chi_{[0, {1 \over 2}]}$ and  $\alpha_1(x)= (2x-1)\chi_{[{1 \over 2},1]}$; then 
\begin{align*}
S_0f&=f \circ \alpha_0 + f \circ \alpha_1, & S_0^*f&={1 \over 2}(f \circ \tau_0+ f \circ \tau_1),
\\
S_1f&=f \circ \alpha_0 - f \circ \alpha_1, & S_1^*f&={1 \over 2}(f \circ \tau_0- f \circ \tau_1 ).
\end{align*}
We choose  unit vectors $\phi$, $\psi$ so that $S_1^*\phi=0$ and $S_0^*\psi=0$. By Cuntz's relation (\ref{eqCuntzRelations}) we
get $S_0S_0^*\phi=\phi$ and $S_1S_1^*\psi=\psi$. Thus we get $\left<\,\phi \mid \psi\,\right>=0$. By our construction,
\[
S^*_0\psi=0
\] 
if
\[
\psi \circ \tau_0= - \psi \circ \tau_1,
\]
 i.e.,
\[
\psi\left(\frac{x}{2}\right)=- \psi\left(\frac{x+1}{2}\right),
\]
 or equivalently,
\[
\psi\left(x\right)=-\psi\left(\frac{2x+1}{2}\right)
=-\psi\left(x+\frac{1}{2}\right)\text{\qquad for all }0 \le x \le \frac{1}{2}.
\]
\end{exmp}

\begin{lem}
\label{LEMMA2.4}
 Let $f$ be an element of $L^2(0,1)$. Then the following statements 
are equivalent:
\begin{enumerate}
\item
\label{LEMMA2.4(1)}
$\displaystyle S_0^*f=0$;
\item
\label{LEMMA2.4(2)}
$\displaystyle f(x)=-f\left(\frac{1}{2}+x\right)$ for all $\displaystyle 0 \le x \le \frac{1}{2}$;
\item
\label{LEMMA2.4(3)}
$\displaystyle f(x)= \sum_{n \ge 0} a_n \cos (2\pi(2n+1)x) + b_n \sin (2\pi(2n+1)x)$.
\end{enumerate}
Moreover $S_1S_1^*$ is the projection onto the closed subspace 
\[
\left\{\,f \in L^2(0,1)\biggm| f(x)=-f\left(\frac{1}{2}+x\right),\;0 \le x \le \frac{1}{2}\,\right\}.
\]

 Similarly the following statements are equivalent:
\begin{enumerate}
\item
\label{LEMMA2.4(1bis)}
$\displaystyle S_1^*f=0$;
\item
\label{LEMMA2.4(2bis)}
$\displaystyle f(x)=f\left(\frac{1}{2}+x\right)$ for all $\displaystyle 0 \le x \le \frac{1}{2}$;
\item
\label{LEMMA2.4(3bis)}
$\displaystyle f(x)=\sum_{n \ge 0} a_n \cos (2\pi(2n)x) + b_n \sin (2\pi(2n)x)$.
\end{enumerate}
Moreover $S_0S_0^*$ is the projection onto the closed subspace 
\[
\left\{\,f \in L^2(0,1)\biggm| f(x)=f\left(\frac{1}{2}+x\right),\;0 \le x \le \frac{1}{2}\,\right\}.
\]
\end{lem}

\begin{pf}
Equivalence of statements (\ref{LEMMA2.4(1)}) and (\ref{LEMMA2.4(2)}) is obvious. That (\ref{LEMMA2.4(3)}) implies (\ref{LEMMA2.4(2)}) is routine as
\[
\cos (x+(2m+1)\pi)=-\cos (x)
\]
and
\[
\sin (x+(2m+1)\pi)=-\sin (x)
\]
for any $m \ge 0$. We will prove now that (\ref{LEMMA2.4(1)}) implies (\ref{LEMMA2.4(3)}). 
To that end, for any $n \ge 0$ we set
\[
c_n(f)= \int^1_0\cos (2\pi nx)f(x)\,dx
\]
and 
\[
s_n(f)=\int^1_0\sin (2\pi nx)f(x)\,dx.
\]
A simple computation shows that $c_n(S_0^*f)=c_{2n}(f)$ and 
$s_n(S_0^*f)=s_{2n}(f)$ for all $n \ge 0$. Thus $c_{2n}(f)=s_{2n}(f)=0$ if $S_0^*f=0$. The last statement 
follows since by the Cuntz relations, $S_0^*f=0$ if and only if $S_1S_1^*f=f$.

For equivalence of statements in the second set, we note that by the Cuntz relations,
$\left<\,f \mid g\,\right>=0$ whenever $S_0^*f=0$ and $S_1^*g=0$. Thus, the
equivalence of the second set of statements follows from that of the first set of
statements. \qed
\end{pf}

\begin{notation}
\label{NotationSubspaces}
We shall need the subspaces
\[
\clk_0:=\left\{\,\psi\mid  S_0^*\psi=0 \,\right\}
\]
and
\[
\left\{\, S_{\Imultiindex}\psi\mid  \Sigma(\Imultiindex ) \le 1 \,\right\}.
\]
The symbol $\Sigma(\Imultiindex )$ in
the second subspace $\{\, S_{\Imultiindex}\psi\mid  \Sigma(\Imultiindex ) \le 1 \,\}$, $\Sigma(I)=\sum_{k}i_{k}$, is defined as all
multi-indices $\Imultiindex$ such that $\Sigma(I) \le 1$. It is motivated as 
follows: In building
ONBs, the aim is to start with a conveniently chosen function $\psi$, and then
to construct the rest from recursively applying monomials in the generators
$S_i$ (chosen from a particular $\clo_N$ representation). This notation allows us
to keep track of the combined system of relations in step-size of length
one.
\end{notation}

For any $\psi$ with $S_0^*\psi=0$, we have
\begin{align*}
\int^1_{\frac12}\psi(x)\psi(2x-1)\,dx 
&= \int^{\frac12}_0 \psi\left(\frac{1}{2} +x\right)\psi(2x)\,dx \\ 
&= - \int^{\frac12}_0\psi(x)\psi(2x)\,dx. 
\end{align*}
Thus we have
\begin{align*}
\left<\,\psi \mid S_1 \psi \,\right> 
&= \int^{\frac12}_0 \psi(x)\psi(2x)\,dx - \int^1_{\frac12}\psi(x)\psi(2x-1)\,dx\\
&= 2 \int^{\frac12}_0\psi(x)\psi(2x)\,dx. 
\end{align*}
More generally for any two elements 
$\psi^{(1)}$, $\psi^{(2)}$
of the subspace
\[
\clk_0:=\left\{\,\psi\mid  S_0^*\psi=0 \,\right\},
\]
we have
\[
\left<\,\psi^{(1)} \bigm| S_1\psi^{(2)}\,\right>=2 \int^{\frac12}_0 
\psi^{(1)}(x)\psi^{(2)}(2x)\,dx.
\] 

Any vector $\psi \in \clh$ such that $S_0^*\psi=0$ and $\left<\,\psi \mid S_1S_0^m\psi\,\right>=0$ for all $m \ge 0$
is called a {\it generating vector} for the closed linear span of the vectors 
$\{\, S_{\Imultiindex}\psi\mid  \Sigma(\Imultiindex ) \le 1 \,\}$
(see Notation \ref{NotationSubspaces} above and Lemma \ref{LEMMA2.6} below for the notation used here). 
A subspace $\clk$ of $\clk_0$ is called 
a {\it basis space} 
for generating vectors if it is a maximal family of vectors that satisfies the following 
mutual relation:
$$\left<\,\psi \mid S_1S^m_0 \psi'\,\right>=0,\qquad m \ge 0,$$ 
for all $\psi,\psi' \in \clk$. Existence of a such a maximal family of vectors follows by Zorn's lemma. 
It is simple to note that $\clk$ is a subspace of $\clk_0$.  

\begin{defn}
\label{Definition2.pound} A maximal subspace $\clk$ as above will be called a
 \emph{basis space}. 
\end{defn}

We now turn to a concrete representation of the subspaces $\clk_0$ and $\clk$ in
$L^2(0,1)$ which were outlined above.

      Our identification of subspaces $\clk_0$ and associated orthonormal bases
is by a certain algorithmic procedure. Below we illustrate our recursive construction
in one 
particular example: In our construction we begin with the representation of
$\clo_2$ from Lemma \ref{LEMMA2.2}, and we give a natural and orthogonal subspace
decomposition of the Hilbert space $L^2(0,1)$. Using this, we then show how an
associated recursive basis may be realized. In our example, we start with a
family of sine functions, normalized to have period one. We then aim for an
orthonormal basis (ONB) when the Hilbert space $L^2(0,1)$ is defined from the
restriction of Lebesgue measure to the unit interval $\Iinterval  = [0, 1]$. (Other
self-similar measures will be considered later!)  Our example will further
serve to illustrate the reflection operations which we will encounter later
in a more general context of self-similar systems.

         As outlined before, the idea is to start our recursion from two
prescribed functions $\psi_0$ and $\psi_1$. Here we take $\psi_0$ to be the constant
function ``one'' on $\Iinterval =[0, 1]$; and we choose $\psi_1 (x) : = s(x) := \sin (2 \pi x)$. The
recursion will be as outlined above: 
The idea is to recursively determine the pair of functions  $(\psi_{2n} ,
\psi_{2n + 1} )$ from sampled subdivisions of $\psi_n$ for each $n \ge 1$. (Note that
our recursion does not begin with $n = 0 $.)  Notation: Set $s_n(x) := \sin( 2
\pi n x) $, $n = 1, 2, \dots$.

\begin{lem}
\label{LEMMA2.5} For every odd integer $n$, the function $s_n$ is in $\clk_0$. Moreover,
if $n$ is even, then $S_0^*s_n$ is non-zero in $L^2(0,1)$.
\end{lem}

\begin{pf}
We begin by setting
$s(x)=\sin (2 \pi x)$, and for any integer $n \ge 1$ we also
set
$
s_n(x)=s(nx)=\sin (2\pi nx)$, $ x \in [0,1]$.
 As
\begin{align*}
s_1\left(\frac{1 + x}{2}\right)&= s\left(\frac{1+x}{2}\right) = 
\sin (\pi(1+x))\\&=- \sin (\pi x)=-s_1\left(\frac{x}{2}\right),
\end{align*}
 we have $s_1 \in \clk_0$. For any odd integer, i.e.,
$n=2m+1$, we check that
\begin{align*}
s_n\left(\frac{1+x}{2}\right) &= \sin (\pi(2m+1)(1+x))=\sin (\pi(2m+1)x + 
(2m+1)\pi)\\&=-\sin (\pi(2m+1)x) = -s_n\left(\frac{x}{2}\right).
\qed
\end{align*}
\end{pf}

We shall need the following additional facts about the functions $s_n$. 
For any two integers $m,n \ge 0 $, we have 
\begin{align*}
\left<\,s_m \mid S_1s_n\,\right> &= 2 \int^{\frac12}_0 s_m(x)s_{n}(2x)\,dx\\ 
&=\int^{\frac12}_0[\cos (2\pi(m-2n)x)-\cos (2\pi(m+2n)x)]\,dx\\
&=\left[\frac{\sin (2\pi(m-2n)x)}{2\pi(m-2n)} - 
\frac{\sin (2\pi(m+2n)x)}{2\pi(m+2n)}\right]{\biggr|}^{\frac12}_0=0 
\end{align*}
 for $m \ne 2n$.
 For $m=2n$, we also check that the integral is $0$.
Since $S^m_0s_n=s_{2^mn}$, we also get 
$\left<\,s_m \mid S_1S_0^m s_n\,\right>=0$. Hence, $\{\, s_{2n+1}(x)\mid  n \ge 0 \,\}$ is a family of orthonormal 
vectors in a basis space. One natural question that we face now: Is it a maximal family, i.e., 
is it a basis space?
We answer this in the
affirmative in the remaining part of this section.

\begin{lem}
\label{LEMMA2.6} Let $\psi$ be a unit vector such that
\[
S_0^*\psi=0\text{\quad and\quad}\left<\,\psi \mid S_1(S_0)^m\psi\,\right>=0
\text{\qquad
for all }m \ge 0,
\]
but
\[
\left<\,\psi \bigm| S_1(S_0)^{m'}S_1 (S_0)^m\psi\,\right> \ne 0\text{\qquad for all }m,m' \ne 0.
\]
Then it follows that the family 
of vectors
\[
\{\, S_{\Imultiindex}\psi\mid \Sigma(\Imultiindex ) \le 1 \,\} ,
\]
where
\[
\Sigma(\Imultiindex ) = \sum_k i_k ,
\]
 and
\[
\Imultiindex =(i_1\,i_2\,\dots\,i_k),\qquad
 i_{k} \in \{ 0,1 \},
\]
is a maximal family of orthonormal vectors. The closed subspace $\clh(\psi)$ spanned by the vectors
$\{\,S_{\Imultiindex}\psi\mid  |\Imultiindex | < \infty \,\}$ 
is invariant under both of the operators $S_0$ and $S_0^*$.   
\end{lem}

\begin{pf}
 Since $S_0^*\psi=0$, it is simple to verify by the Cuntz relations that $(S_0)^m\psi$ is 
orthogonal to $(S_0)^n\psi$ for $m \ne n \ge 0$, where by convention $(S_0)^0=\Iidentity $. As $S_1^*S_0=0$ and 
$S_0^*\psi=0$, we also check that $(S_0)^{m'}S_1(S_0)^m\psi$ is orthogonal to $(S_0)^{n'}S_1(S_0)^{n}\psi$ 
and $(S_0)^{n'}\psi$ for all $m' \ne n'$ and $n,n' \ge 0$. Thus we are left to check for $m'=n'$. In such 
case orthogonality follows by our hypothesis that $\left<\,\psi \mid S_1(S_0)^m\psi\,\right>=0$ for all $m \ge 0$. It is clear 
that the vector space generated by these vectors is both $S_0$- and $S_0^*$-invariant. The maximal property 
is also evident. See Section \ref{SNew3} below for details. \qed
\end{pf}

Terminology: $\mathcal{M}$ will denote the set of all finite multi-indices.

\begin{prop}
\label{PROPOSITION2.7} Let $(S_i)$ be the irreducible representation of $\clo_2$ as in Example \textup{\ref{EXAMPLE2.3}}, and let
$\clk$ be a basis space. 
For a vector $\psi$ use the notation $\clh(\psi)$ as in Lemma \textup{\ref{LEMMA2.6}} for the closed subspace spanned by the vectors
$\{\,S_{\Imultiindex}\psi\mid  |\Imultiindex | < \infty \,\}$.
Then the following hold:
\begin{enumerate}
\item 
\label{PROPOSITION2.7(1)}
For each unit vector $\psi \in \clk$,
the vectors in the family
$\{\, S_{\Imultiindex}\psi\mid \Sigma (\Imultiindex ) \le 1 \,\}$ are
orthonormal, i.e., norm one, and
mutually orthogonal.  
\item 
\label{PROPOSITION2.7(2)}
For any two orthogonal unit vectors $\psi^{(1)},\psi^{(2)} \in \clk$, 
$\clh(\psi^{(1)})$ is orthogonal to $\clh(\psi^{(2)})$. 
\end{enumerate}
\end{prop}

\begin{pf}
Proof is routine as in Lemma \ref{LEMMA2.6}. \qed
\end{pf}

\begin{lem}
\label{LEMMA2.8}
 \textup{(}Haar--Walsh, Fig.\ \textup{\ref{FigHaarSequence}} above, \textup{\cite{Haa10,Wal23})}
 Let $\clh=L^2(0,1)$ and let $\phi_0=\chi_{[0,1]}$. Then the recursive system
\begin{align*}
\phi_{2n}(x)&= \phi_n(2x)+\phi_n(2x-1),
\\
\phi_{2n+1}(x)&=\phi_n(2x)-\phi_n(2x-1)
\end{align*}
for $n \ge 0$ defines an orthonormal basis for $\clh$.
\end{lem}

\begin{pf}
 We set $S_0,S_1$ as in Example \ref{EXAMPLE2.3}. We consider the Hardy space $\clh_{+}$ given by
$$f(z)=\sum_{n \ge 0} c_n z^n,\qquad z \in \IT^1= \{\,z \in \IC\mid |z|=1 \,\},$$
for $(c_n) \in l^2$ and $\left\Vert f\right\Vert ^2=\sum_{n \ge 0}|c_n|^2$. We also set $e_n(z)=z^n$ and 
\begin{align*}
\tilde{S}_0f(z)&=f(z^2),
\\
\tilde{S}_1f(z)&=zf(z^2)
\end{align*}
for all $f \in \clh_{+}$ and $z \in \IT^1$. We have $S^*_0\phi_0=\phi_0$ and $\tilde{S}^*_0e_0=e_0$. Now, defining
$W\colon\clh_{+} \raro \clh$ by $We_n=\phi_n$, we verify that
\[
W\tilde{S}_i=S_iW\text{\qquad for }i\in\{0,1\}.
\]
If $n=j_1+j_22+\dots+j_k2^k$ where $j_r \in \{0,1\}$ is the dyadic representation of an integer 
$n \ge 0$, it follows that $We_n=W\tilde{S}_{j_1}\tilde{S}_{j_2}\cdots\tilde{S}_{j_k}e_0=S_{j_1}S_{j_2}\cdots S_{j_k}\phi_0=\phi_n$,
and the result follows.
In the last step we are using the fact (details in \cite{Jor05})
that the two $\clo_2$ representations are irreducible; so the intertwining
operator $W$ is a constant times a unitary.
 \qed
\end{pf}


\begin{definitions}
\label{Def}\setcounter{enumd}{0}\refstepcounter{enumd}\label{Def.a}(\ref{Def.a})~For
$\psi\in\mathcal{H}\setminus\left\{  0\right\}  $, set%
\begin{multline*}
K\left(  \psi\right)  :=\min\{\,K\in\mathcal{M}\mid\mathrm{s.t.}\left\langle
\,S_{J_{1}}\psi\mid S_{J_{2}}\psi\,\right\rangle =0\;\forall\,J_{i}\lvertneqq
K,\;J_{1}\neq J_{2}\text{,}\\
\text{ and }\exists\,J\lvertneqq K\;\mathrm{s.t.}\;\left\langle \,S_{J}%
\psi\mid S_{K}\psi\,\right\rangle \neq0\,\}.
\end{multline*}

\refstepcounter{enumd}\label{Def.b}(\ref{Def.b})~Further, define%
\begin{equation}
\mathcal{H}\left(  \psi\right)  :=\operatorname{closed\;span}\left\{  \,S_{0}^{m}%
S_{J}\psi\mid m\in\mathbb{N}_{0},\;J\lvertneqq K\left(  \psi\right)
\,\right\}  . \label{eqDef.b}%
\end{equation}

\refstepcounter{enumd}\label{Def.c}(\ref{Def.c})~A family of vectors $\left(
\psi_{n}\right)  $, $\left\Vert \psi_{n}\right\Vert =1$, is said to be
\emph{maximal} and orthogonal if
\setcounter{enumdd}{0}\refstepcounter{enumdd}\label{Def.ci}(\ref{Def.ci})~the
corresponding subspaces $\mathcal{H}\left(  \psi_{n}\right)  $ are mutually
orthogonal, \refstepcounter{enumdd}\label{Def.cii}(\ref{Def.cii}%
)~$\mathcal{H}\left(  \psi_{n}\right)  \perp\varphi$, and
\refstepcounter{enumdd}\label{Def.ciii}(\ref{Def.ciii})~they are not part of a
bigger such family.
\end{definitions}

\begin{thm}
\label{THEOREM2.8} 
Let $\clk$ be a basis space, and let $\{\, \psi_n \mid n \ge 1 \,\}$ be an orthonormal basis
for $\clk$. Further, let $\clh_n$ be a maximal family of subspaces associated with
$(\psi_n)$
in
the sense of Definition \textup{\ref{Def}(\ref{Def.c})}.
Then 
\[
\clh= \IC \, \phi \bigoplus_{n \ge 1} \clh_n.
\]
\end{thm}

Before turning to the proof, we will need some preliminaries which are
included in the remarks below. The proof will then be resumed at the end of 
Section \ref{SNew3}
below.

       The reasoning in the proof of the theorem has the following two parts
in rough outline: Firstly, the argument for orthogonality of the spaces and
the vectors which go into our basis construction is largely combinatorial,
and it is sketched in Remarks \ref{REMARKS2.9} below. 

       Secondly, we must prove that the functions which are produced by the
algorithm form a total family in $L^2(0,1)$, i.e., that these vectors span a
dense subspace in $L^2(0,1)$. Recall our algorithm starts with the two
functions $\psi_0 = \phi = {}$the constant function ``one'', and $\psi_1(x) = s_1(x) =
\sin ( 2 \pi x)$. Using the sine functions from Lemma \ref{LEMMA2.5}, and our
representation of $\clo_2$ from Lemmas \ref{LEMMA2.2} and \ref{LEMMA2.6}, we then organize a system of
orthogonal functions in each of the closed subspaces $\clh_n$ from the conclusion
of the theorem. Our assertion is that this orthonormal family of vectors is
total. Our argument going into the proof of this is structured as follows
(details in Section \ref{SNew3}): Suppose some $f$ in $L^2(0,1)$ is in the orthogonal
complement of the family. We then show that both vectors $f$ and $S_1^* f$ must
have Fourier expansions consisting only of cosine functions. Because of the
reflection built into the operator $S_1^* $, we conclude that this is only
possible if $f$ is zero.

\begin{remarks}
\label{REMARKS2.9} Our algorithm may be applied both to existing wavelets, and to new ones as well. Consider for example the Haar--%
Walsh sequence
\[
\{\, \phi_n\mid n \ge 0 \,\} \subseteq L^2(0,1)
\]
defined as in Lemma \textup{\ref{LEMMA2.8}}. For any multi-indices $J=(j_1\,j_2\,\dots\,j_n) \in \clm$,
$j_i \in \{0,1\}$, consider the orthogonal family $S_J\,\phi$ and set
$$\clm_{ev}=\{\, J \in \clm\mid \Sigma(J)=0,2,4,6, \dots \,\}$$
and 
\be
\clm_1= \{\, J \in \clm\mid \Sigma(J) \le 1 \,\}.
\label{eq2.4}
\ee
    Here a number $N$ is fixed and our multi-indices are built from the
alphabet $\{0, 1, \dots, N -1 \}$. Further we define $\Sigma(J)$ to be $\sum_1^n j_i  $.

For any $J \in \clm_{ev}$, we set
$$\clh(S_J\psi)=\operatorname{closed\;span}\{\,S_{KJ}\psi\mid K \in \clm_1 \,\}.$$
We claim that
\be
\sideset{}{^{\smash{\oplus}}}{\sum}\limits
_{J \in \clm_{ev}} \clh(S_J\psi) \oplus \IC \, \phi_0 = L^2(0,1)
\label{eq2.5}
\ee
and the vectors appearing in (\ref{eq2.5}) are orthonormal. Note that the first few terms in the system of closed subspaces in (\ref{eq2.5})
are
$$\clh(\phi_1) \oplus \clh(\phi_7) \oplus \clh(\phi_{11}) \oplus \clh(\phi_{13}) \oplus\cdots.$$
Specifically,
\begin{align*}
\phi_7&=S_1^2\psi,\\
\phi_{11}&=S_1^2S_0\psi,\\
\phi_{13}&=S_1S_0S_1\psi,\\
\vdots&\,.
\end{align*} 
\end{remarks}

              We now turn to our basis constructions.

              In brief outline (Proof of Theorem \ref{THEOREM2.8} Part 1): The 
starting point is a fixed IFS with $N$
endomorphisms, and an associated representation $(S_i)$ of $\clo_N$ in
$L^2(\mu)$. Although it would seem natural to begin with cyclic subspaces in
$L^2(\mu)$,  care must be exercised in order to guarantee orthogonality. Our
starting point will be two orthogonal (and carefully selected) vectors $\phi$
and $\psi$ satisfying $S_0 \phi = \phi$, and $S_0^* \psi = 0$. (Note that
orthogonality of $\phi$ and $\psi$ is implied by these relations.) Then using the
representation, new vectors $\phi$ are constructed to make up part of an ONB.
The new vectors result from applying operator monomials in the generators
$S_i$  to $\psi$.  So operator monomials are applied to the single vector $\psi$,
and the monomials are selected with view to orthogonality. The various new
orthogonal vectors $\phi$ are assigned subscripts according to the particular
operator monomial which is applied to $\psi$. The process is then repeated
inductively on additional vectors $\psi$ as needed for creating a full ONB.

\newsection{\label{SNew3}Irreducible representations of $\clo_N$}

Let $N\in\mathbb{N}$, $N\geq2$, and let $\mathcal{H}$ be a separable Hilbert
space. Let $\left\{  \,S_{i}\mid0\leq i<N\,\right\}  $ be a system of
isometries in $\mathcal{H}$ which define an irreducible representation of the
Cuntz algebra $\mathcal{O}_{N}$. We shall further assume that there is some
$\varphi\in\mathcal{H}$ such that $\left\Vert \varphi\right\Vert =1$ and
$S_{0}\,\varphi=\varphi$. (Note that then $S_{0}^{\ast}\,\varphi=\varphi$ as well.)

Set $\mathcal{M=M}_{N}$ ($={}$the set of all finite multi-indices $J=\left(
j_{0}\,j_{1}\,\dots\,j_{p}\right)  $ where $j_{i}\in\left\{  0,1,\dots
,N-1\right\}  $). Since%
\[
\left(  j_{0}\,j_{1}\,\dots \,j_{p}\right)  \longmapsto j_{0}+j_{1}N+\dots
+j_{p}N^{p}%
\]
defines a bijection of $\mathcal{M}_{N}$ onto $\mathbb{N}_{0}=\left\{
0,1,\dots\right\}  $, $\mathcal{M}_{N}$ acquires an order induced from the
natural order on $\mathbb{N}_{0}$.

\begin{thm}
\label{ThmNew3.1}Let $N$ be an integer${} \ge 2$, let $\mathcal{H}$
be a Hilbert space, and let $\varphi
\in\mathcal{H}$, $\left\Vert \varphi\right\Vert =1$, be chosen such that
$S_{0}\varphi=\varphi$, where $\left(  S_{i}\right)  _{i=0}^{N-1}$ is a given
irreducible representation of $\mathcal{O}_{N}$ on $\mathcal{H}$. Then there
is a maximal family $\left(  \psi_{n}\right)  _{n\geq1}$, such that
$\left\Vert \psi_{n}\right\Vert =1$, $S_{0}^{\ast}\psi_{n}=0$; and%
\begin{equation}
\mathcal{H}=\mathbb{C}\,\varphi\oplus
\sideset{}{^{\smash{\oplus}}}{\sum}\limits
_{n\geq1}\mathcal{H}\left(
\psi_{n}\right)  .\label{eqThmNew3.1}%
\end{equation}
Note, every maximal family will satisfy \textup{(\ref{eqThmNew3.1})}.
\end{thm}

Before starting the proof, we need a lemma.

\begin{lem}
\label{LemNew3.2}
Let $\psi\in\mathcal{H}\setminus
\left\{  0\right\}  $. Assumptions as above, including $S_0^* \psi = 0$.
Then the corresponding subspace $\mathcal{H}
\left(  \psi\right)  $ is invariant under both of the
operators $S_{0}$ and $S_{0}^{\ast}$.
\end{lem}

\begin{pf}
By
\textup{(\ref{eqDef.b})}, $\mathcal{H}\left(  \psi\right)  $ is the closed
span of $\left\{  \,S_{0}^{m}S_{J}\psi\mid m\in\mathbb{N}_{0},\;J\lvertneqq
K\left(  \psi\right)  \,\right\}  $. Since $S_{0}^{\ast}\psi=0$, the only
non-trivial case to consider is the set of vectors $S_{0}^{\ast}S_{J}\psi$
when $J\neq\varnothing$ and $J\lvertneqq K\left(  \psi\right)  $. But note
that $S_{0}^{\ast}S_{J}\psi=S_{J^{\prime}}\psi$ where $J^{\prime}\lvertneqq
K\left(  \psi\right)  $. This concludes the proof.\qed
\end{pf}

\begin{pf*}{PROOF OF THEOREM \textup{\ref{ThmNew3.1}}}
Using Zorn's lemma, we may select a
maximal family $\left(  \psi_{n}\right)  $ as stated in the theorem. When it
is chosen, our claim is that then (\ref{eqThmNew3.1}) holds.

We finish the argument by assuming that the orthocomplement $\mathcal{L}$
of $\mathbb{C}\,\varphi\oplus\sum_{n\geq1}\mathcal{H}\left(  \psi_{n}\right)  $,
$\mathcal{L} =
\left\{ \mathbb{C}\,\varphi\oplus\sum_{n\geq1}\mathcal{H}\left(  \psi_{n}\right)  \right\}_{\phantom{\perp}}^{\perp} = \mathcal{H} \ominus \left\{ \mathbb{C}\,\varphi\oplus\sum_{n\geq1}\mathcal{H}\left(  \psi_{n}\right)  \right\}$,
is nonzero, and then derive a contradiction.

By Lemma \ref{LEMMA2.2}, there are two cases for the projection $S_{0}%
S_{0}^{\ast}$:  

\textbf{Case 1.} $S_{0}S_{0}^{\ast}|_{\mathcal{L}}\lvertneqq \Iidentity _{\mathcal{L}}$; and

\textbf{Case 2.} $S_{0}S_{0}^{\ast}|_{\mathcal{L}}=\Iidentity _{\mathcal{L}}$.

Assuming Case 1, there is some $\psi\in\mathcal{L}\setminus\left\{  0\right\}
$ such that $S_{0}^{\ast}\psi=0$. Following the argument from 
the proof of
Theorem \ref{THEOREM2.8}, we now
prove that%
\[
\mathcal{H}\left(  \psi\right)  \perp\mathcal{H}\left(  \psi_{n}\right)
\text{\qquad for all }n\geq1;
\]
and this then contradicts maximality of the family $\left(  \psi_{n}\right)  $.

Assuming Case 2, and $\mathcal{L}\neq0$, we get%
\[
S_{i}S_{i}^{\ast}\psi=0\text{\qquad for all }i\text{, }1\leq i<N-1\text{, and
all }\psi\in\mathcal{L}.
\]
This follows from%
\[
\sum_{i=0}^{N-1}S_{i}S_{i}^{\ast}\psi=\psi.
\]
We conclude that $S_{i}^{\ast}\psi=0$ for $i=1,2,\dots,N-1$.

Now use the properties of the vector $\varphi$ to conclude that%
\begin{equation}
\left\langle \,S_{J}^{\ast}\psi\mid\phi\,\right\rangle =\left\langle
\,\psi\mid S_{J}\varphi\,\right\rangle =0\text{\qquad for all }J\in
\mathcal{M}.\label{eqThmNew3.1p}%
\end{equation}
We have $S_{0}\varphi=S_{0}^{\ast}\varphi=\varphi$, and therefore $S_{i}%
^{\ast}\varphi=0$ for $i\geq1$. As a result, the closed subspace spanned by
$\left\{  \,S_{J}\varphi\mid J\in\mathcal{M}_{N}\,\right\}  $ is invariant
under $\mathcal{O}_{N}$. Since $\varphi\neq0$, this subspace must be all of
$\mathcal{H}$, and we get $\psi=0$ from (\ref{eqThmNew3.1p}), which is a contradiction.

This concludes the proof of the theorem.
\qed
\end{pf*}

\textbf{Proof of Theorem \ref{THEOREM2.8}} concluded.  We have chosen to 
introduce the tools of Cuntz algebra representations before the conclusion of 
the proof of Theorem \ref{THEOREM2.8}. Once this is available, the reader will 
easily be able to establish the existence of a maximal family of subspaces as 
in the conclusion of Theorem \ref{THEOREM2.8}, following the reasoning above 
(in the proof of Theorem \ref{ThmNew3.1}) using Zorn, and in the same way 
dividing the analysis up into two cases.

\newsection{\label{SNew4}A scale-$4$ Cantor set}

In this section we continue work on harmonic analysis for iterated-function-%
system fractals started in \cite{JoPe98} and continued in \cite{DuJo03}. Our present
basis constructions differ from the earlier ones mainly in their emphasis on
the use of ``localized'' functions.

\begin{figure}[pt]
\begin{center}
\includegraphics[bb=0 86 288 202,height=116bp,width=288bp]{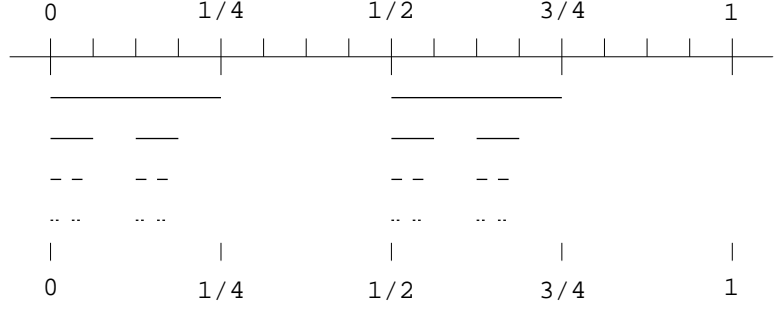}
\end{center}
\caption{The Cantor set $X_{4}$}%
\label{FigCantor}%
\end{figure}

        We begin the section with an outline of a new application of our
general algorithm from
Theorem \ref{ThmNew3.1} to the case of the scale-$4$ Cantor set $X=X_{4}$
on the real line $\mathbb{R}$. This Cantor set was studied in \cite{JoPe98}. We
begin with a brief review of $X_{4}$ and its associated measure $\mu=\mu_{4}$:
First note that by \cite{Hut81}, or by direct computation, there is a unique Borel
probability measure $\mu$ on $\mathbb{R}$ satisfying
\begin{equation}
\int f\,d\mu=\frac{1}{2}\left(  \int f\left(  \frac{x}{4}\right)
\,d\mu\left(  x\right)  +\int f\left(  \frac{x+2}{4}\right)  \,d\mu\left(
x\right)  \right)  \label{eqNew4.1}%
\end{equation}
for all bounded continuous functions $f$. Its support $X$ is the unique
compact subset of $\mathbb{R}$ satisfying $X=\tau_{0}\left(  X\right)
\cup\tau_{1}\left(  X\right)  $ where
\begin{equation}
\tau_{0}\left(  x\right)  =\frac{x}{4}\,,\qquad\tau_{1}\left(  x\right)
=\frac{x+2}{4}\,.
\label{eqUniqueCompact}
\end{equation}
The set $X$ is sketched in Fig.\ \ref{FigCantor}.
Specifically, $X = \left\{\,
\sum_{i=1}^{\infty} k_i/4^i  \mid k_i \in \{0, 2\} \,\right\}$. 
When fractions inside $[0, 1]$ are
written in base $4$, $X$ results as a Cantor set from always omitting the two
choices $1$ and $3$. Equivalently (see Fig.\ \ref{FigCantor}), by repeated subdivision of $[0,
1]$, our Cantor set $X$ results when, at each quarter subdivision step, we are
omitting the second and the fourth open subintervals. 

It is evident that $X\subset\left[  0,1\right]  $ and that
$\operatorname*{supp}\left(  \mu\right)  =X$ where $\operatorname*{supp}%
\left(  \mu\right)  $ denotes the support of $\mu$. Furthermore, the measure
$\mu$ is the restriction of the Hausdorff measure $\operatorname*{Haus}%
_{\frac{1}{2}}$ of Hausdorff dimension $\frac{1}{2}$. We shall be interested
in the algorithms for ONBs in the Hilbert space $L^{2}\left(  \mu\right)
=L^{2}\left(  X_{4},\mu_{4}\right)  $. One reason for the choice of this
particular example is that it appears to be prototypical of IFSs for which
$L^{2}\left(  \mu\right)  $ has an ONB consisting of complex exponentials
$e_{\lambda}\left(  x\right)  :=e^{i2\pi\lambda x}$. Specifically, the
co-authors of \cite{JoPe98} prove that $\left\{  \,e_{\lambda}\mid\lambda
\in\left\{  0,1,4,5,16,17,20,21,\dots\right\}  \,\right\}  $ is an ONB in
$L^{2}\left(  \mu\right)  $. (Note that the index set for this ONB is the set
of all finite sums%
\[
\left\{  \,k_{0}+k_{1}4+k_{2}4^{2}+\dots+k_{s}4^{s}\bigm| k_{i}\in\left\{
0,1\right\}  \,\right\}  \text{.)}
\]
There are other choices of ONBs, see, e.g., \cite{DuJo05b} and \cite{DuJo05c}, but the argument for
orthogonality is based on a form of factorization of the following transform:%
\begin{align*}
\hat{\mu}\left(  \lambda\right)   &  :=\int e_{\lambda}\left(  x\right)
\,d\mu\left(  x\right) \\
&  =\frac{1}{2}\left(  1+e^{i\pi\lambda}\right)  \hat{\mu}\left(
\frac{\lambda}{4}\right) \\
&  =\prod_{m=0}^{\infty}\frac{1}{2}\left(  1+e^{i\pi\lambda4^{-m}}\right)
,\qquad\lambda\in\mathbb{R}.
\end{align*}
The two operators $S_{0}$, $S_{1}$ generating our associated representation of
$\mathcal{O}_{2}$ in $L^{2}\left(  X_{4},\mu_{4}\right)  $ are specified as
follows:%
\begin{equation}
\left(  S_{0}f\right)  \left(  x\right)  =\left\{
\begin{aligned}&f\left( 4x\right) ,&&\qquad x\in\tau_{0}\left( X_{4}\right) ,\\&f\left( 4x-2\right) ,&&\qquad x\in\tau_{1}\left( X_{4}\right) ,\end{aligned}\right.
\label{eqNew4.S0}%
\end{equation}
and%
\begin{equation}
\left(  S_{1}f\right)  \left(  x\right)  =\left\{
\begin{aligned}&f\left( 4x\right) ,&&\qquad x\in\tau_{0}\left( X_{4}\right) ,\\&-f\left( 4x-2\right) ,&&\qquad x\in\tau_{1}\left( X_{4}\right) .\end{aligned}\right.
\label{eqNew4.S1}%
\end{equation}
We have $S_{0}e_{0}=e_{0}$, and the formula (\ref{eqNew4.1}) takes the form%
\[
\int_{X}f\,d\mu=\int_{X}\left(  S_{0}^{\ast}f\right)  \left(  x\right)
\,d\mu\left(  x\right)  ,\qquad f\in C\left(  X\right)  .
\]

We now prove that the representation of $\mathcal{O}_{2}$ outlined above for
the Cantor example is \emph{irreducible}.

\begin{thm}
\label{ThmNew.4p}Let $\mathcal{H}=L^{2}\left(  X_{4},\mu_{4}\right)  $ be the
Hilbert space defined from the Cantor measure $\mu_{4}$ with
$\operatorname*{supp}\left(  \mu_{4}\right)  =X_{4}$. Let $S_{i}$, $i=0,1$, be
the operators defined in \textup{(\ref{eqNew4.S0})} and
\textup{(\ref{eqNew4.S1})} above.

These two operators then generate an irreducible representation of the Cuntz
$C^{\ast}$-algebra $\mathcal{O}_{2}$ acting on the Hilbert space $\mathcal{H}$.
\end{thm}

\begin{pf}
We first check that the two operators $S_{i}$, $i=0,1$, are isometries. The
calculation uses the three facts:

\begin{enumerate}
\renewcommand{\theenumi}{\roman{enumi}}

\item \label{ThmNew.4p(1)}$\tau_{0}\left(  X\right)  \cap\tau_{1}\left(
X\right)  =\varnothing$,

\item \label{ThmNew.4p(2)}$\tau_{0}\left(  X\right)  \cup\tau_{1}\left(
X\right)  =X$,

\item \label{ThmNew.4p(3)}$\mu\left(  cE\right)  =\sqrt{c}\,\mu\left(
E\right)  $ for all $c\in\mathbb{R}_{+}$ and all $E\in\mathcal{B}_{\frac{1}%
{2}}\left(  X\right)  ={}$the sigma-algebra of the $\operatorname*{Haus}%
_{\frac{1}{2}}$-measurable sets.
\end{enumerate}

Now if $f$ is bounded and measurable, then%
\begin{align*}
\int_{X}\left\vert S_{i}f\right\vert ^{2}\,d\mu &  =\int_{\tau_{0}\left(
X\right)  }\left\vert f\left(  4x\right)  \right\vert ^{2}\,d\mu\left(
x\right)  +\int_{\tau_{1}\left(  X\right)  }\left\vert f\left(  4x-2\right)
\right\vert ^{2}\,d\mu\left(  x\right)  \\
&  =\int_{X}\left\vert f\left(  y\right)  \right\vert ^{2}\,d\mu\left(
\tau_{0}\left(  y\right)  \right)  +\int_{X}\left\vert f\left(  y\right)
\right\vert ^{2}\,d\mu\left(  \tau_{1}\left(  y\right)  \right)  \\
&  =\frac{1}{2}\int_{X}\left\vert f\right\vert ^{2}\,d\mu+\frac{1}{2}\int
_{X}\left\vert f\right\vert ^{2}\,d\mu=\int_{X}\left\vert f\right\vert
^{2}\,d\mu.
\end{align*}

The remaining Cuntz property, $S_{0}S_{0}^{\ast}+S_{1}S_{1}^{\ast
}=\Iidentity _{\mathcal{H}}$, is immediate from
\[
\mu=\frac{1}{2}\left(  \mu\circ\tau_{0}^{-1}+\mu\circ\tau_{1}^{-1}\right)  .
\]

Set $e_{0}=\varphi=\chi_{X}={}$the indicator function of the Cantor set
$X=X_{4}$. We claim that the vectors $S_{J}\varphi$, $J\in\mathcal{M}$ ($={}%
$all finite $0$--$1$ multi-indices), span a dense subspace in $\mathcal{H}$.
Since $S_{0}^{\ast}\varphi=\varphi$ and $S_{1}^{\ast}\varphi=0$, it follows
that $\varphi$ is then a cyclic vector for the $\mathcal{O}_{2}$%
-representation; and that $\omega_{\varphi}:=\left\langle \,\varphi\mid
\cdot\;\varphi\,\right\rangle $ is a Cuntz state. An application of 
\cite[Theorem 3.3]{BJP96}
 then yields irreducibility.

For all multi-indices $J$, $K$ of the same length, we set $JK:=\sum_{i}%
j_{i}k_{i}$. Recall $j_{i},k_{i}\in\left\{  0,1\right\}  $ for $\mathcal{O}%
_{2}$. An induction now yields%
\begin{equation}
\chi_{\tau_{J}\left(  X\right)  }=\frac{1}{\left\vert J\right\vert }%
\sum_{\substack{K,\\\left\vert K\right\vert =\left\vert J\right\vert }}\left(
-1\right)  ^{JK}S_{K}\varphi.
\label{eqchitauJ}
\end{equation}
But by the Cantor construction, we know that these indicator functions
\begin{equation}
\left\{  \,\chi_{\tau_{J}\left(  X\right)  }\bigm| J\in\mathcal{M}\,\right\}  
\label{eqIndicator}
\end{equation}
span a dense subspace in $\mathcal{H}=L^{2}\left(  \mu_{4}\right)  $.
The fact that the functions in (\ref{eqIndicator}) are
total in $L^{2}\left(  \mu_{4}\right)  $ is a consequence of the recursive algorithm used in the
construction of the measure $\mu_{4}$ (see Fig.\ \ref{FigCantor}). Moreover, the argument for the
present special case in fact carries over \emph{mutatis mutandis} to general IFS
constructions, regardless of whether the limit measure $\mu$ is fractal or not;
see Section \ref{S2} above and \cite{Hut81} for further details.

This concludes the proof of our theorem.
\qed
\end{pf}

\begin{cor}
\label{CorNew4.pp}\textup{\cite{JoPe98}} Consider the Cantor system $\left(
X,\mu\right)  =\left(  X_{4},\mu_{4}\right)  $, the Hilbert space
$\mathcal{H}=L^{2}\left(  X,\mu\right)  $, and the functions
\begin{align}
e_{n}\left(
x\right)  &:=e^{i2\pi nx}, \nonumber \\
n&\in\left\{  0,1,4,5,16,17,20,21,\dots\right\}  \nonumber \\
&=\left\{  \,j_{0}+j_{1}4+j_{2}4^{2}+\dots+j_{p}4^{p}\bigm| j_{i}\in\left\{
0,1\right\}  \,\right\}  =:\Lambda.
\label{eqCorNew4.pp.1}
\end{align}
 These functions form an ONB in
$\mathcal{H}$, and%
\begin{equation}
\mathcal{H}=\mathbb{C}\,\varphi\oplus
\sideset{}{^{\smash{\oplus}}}{\sum}\limits
_{\substack{m\colon \\m\text{ odd,}
\\m\in\Lambda}}\mathcal{H}\left(  e_{m}\right)  .\label{eqCorNew4.pp.2}%
\end{equation}

\end{cor}

\begin{pf}
The result is a direct corollary of the two theorems, Theorem \ref{ThmNew3.1}
and Theorem \ref{ThmNew.4p}. In the present application it turns out that
$K\left(  e_{m}\right)  =0$ for all $m$ odd in $\Lambda$. As a result we get
that $\mathcal{H}\left(  e_{m}\right)  $ is spanned by $S_{0}^{k}e_{m}$,
$k\in\mathbb{N}_{0}$; and recall that $S_{0}^{k}e_{m}=e_{m4^{k}}$, $m$ odd in
$\Lambda$.
\qed
\end{pf}

\begin{ack}
We thank Brian Treadway for excellent typesetting and for
producing the figures, and for suggesting a number of improvements. Finally,
we are grateful to two anonymous referees for offering very helpful
clarifications and suggestions.  And we thank Prof Myung-Sin Song for help 
with revisions.
\end{ack}


\begin{thebibliography}{ALTW04}
\frenchspacing
\bibitem[ALTW04]{ALTW04}
Akram Aldroubi,  David Larson,  Wai-Shing Tang, and Eric Weber, 
{\em Geometric aspects of frame representations of abelian groups}, Trans. Amer.
Math. Soc. {\bf 356} (2004), no. 12, 4767--4786.

\bibitem[Ash90]{Ash90}
Robert B. Ash, {\em Information {T}heory}, Dover, New York, 1990, corrected reprint
  of the original 1965 Interscience/Wiley edition.

\bibitem[BJMP05]{BJMP05}
L. W. Baggett, P. E. T. Jorgensen, K. D. Merrill, and J. A. Packer, {\em
  Construction of {P}arseval wavelets from redundant filter systems}, J. Math.
  Phys. {\bf 46} (2005), no.~8, 083502, 28 pp., doi:10.1063/1.1982768.

\bibitem[BPS03]{BPS03}
Ernst Binz, Sonja Pods, and Walter Schempp, {\em Heisenberg
groups---a unifying structure of signal theory, holography and quantum
information theory}, J. Appl. Math. Comput. {\bf 11} (2003), no. 1--2, 1--57.

\bibitem[BoPa05]{BoPa05}
Bernhard G. Bodmann and  Vern I. Paulsen, {\em Frames, graphs and
erasures}, Linear Algebra Appl. {\bf 404} (2005), 118--146.

\bibitem[BJP96]{BJP96}
Ola Bratteli, Palle E. T. Jorgensen, and Geoffrey L. Price, {\em Endomorphisms of
  $\mathcal{B}(\mathcal{H})$}, Quantization, Nonlinear Partial Differential
  Equations, and Operator Algebra (W.~Arveson, T.~Branson, and I.~Segal, eds.),
  Proc. Sympos. Pure Math., vol.~59, American Mathematical Society, Providence,
  1996, pp.~93--138.

\bibitem[Chr55]{Chr55}
H. E. Chrestenson, {\em A class of generalized Walsh functions},
Pacific J. Math. {\bf 5} (1955), no. 1, 17--31.

\bibitem[CoWi92]{CoWi92}
R. R. Coifman and M. V. Wickerhauser {\em Entropy-based algorithms
for best basis selection}, IEEE Trans. Information Theory {\bf 38} (1992) 713--718.

\bibitem[Cu77]{Cu77}
Joachim Cuntz,
{\em Simple $C\sp*$-algebras generated by isometries},
Comm. Math. Phys. {\bf 57} (1977), no. 2, 173--185.

\bibitem[Dau92]{Dau92}
Ingrid Daubechies, {\em Ten Lectures on Wavelets}, CBMS-NSF Regional
Conference Series in Applied Mathematics, vol. 61, Society for Industrial and
Applied Mathematics (SIAM), Philadelphia, PA, 1992.

\bibitem[Dut05]{Dut05}
Dorin Ervin Dutkay, {\em Some equations relating multiwavelets and
multiscaling functions}, J. Funct. Anal. {\bf 226} (2005), no. 1, 1--20.

\bibitem[DuJo05]{DuJo05b}
Dorin Ervin Dutkay and Palle E. T. Jorgensen, {\em Hilbert spaces of martingales supporting
  certain substitution-dynamical systems}, Conform. Geom. Dyn. {\bf 9} (2005),
  24--45.

\bibitem[DuJo06a]{DuJo05c}
Dorin Ervin Dutkay and Palle E. T. Jorgensen, 
{\em Iterated function systems, Ruelle operators,
 and invariant projective measures}, Math. Comp. {\bf 75} (2006),
  1931--1970.

\bibitem[DuJo06b]{DuJo03}
 Dorin Ervin Dutkay and Palle E. T. Jorgensen, 
{\em Wavelets on fractals}, Rev. Mat.
  Iberoamericana {\bf 22} (2006), no. 1, 131--180.

\bibitem[Gli60]{Gli60}
James G. Glimm, {\em On a certain class of operator algebras}, Trans. Amer. Math.
  Soc. {\bf 95} (1960), 318--340.
 
\bibitem[Haa10]{Haa10}
Alfred Haar, {\em Zur {T}heorie der orthogonalen {F}unktionensysteme}, Math. Ann.
  {\bf 69} (1910), 331--371.

\bibitem[Hut81]{Hut81}
John E. Hutchinson,
{\em Fractals and self similarity}, Indiana Univ. Math.~J.
  {\bf 30} (1981), 713--747.

\bibitem[JMR01]{JMR01}
St{\'e}phane Jaffard, Yves Meyer, and Robert D. Ryan, 
{\em Wavelets: Tools for {S}cience \&
  {T}echnology}, revised ed., SIAM, Philadelphia, 2001.

\bibitem[Jor05]{Jor05}
Palle E. T. Jorgensen, {\em Measures in wavelet decompositions},
Adv. Appl. Math. {\bf 34} (2005), no. 3, 561--590.

\bibitem[Jor06]{Jor06}
Palle E. T. Jorgensen, {\em Analysis and {P}robability: Wavelets, {S}ignals,
  {F}ractals}, Grad. Texts in Math., vol. 234, Springer, New York, 2006.

\bibitem[JoOl98]{JoOl98}
  Palle E. T.  Jorgensen
and
Gestur \'Olafsson,
  {\em Unitary representations of Lie groups with reflection
symmetry},
 J. Funct. Anal. {\bf 158} (1998), no.
  1, 26--88.


\bibitem[JoPe94]{JoPe94}
Palle E. T. Jorgensen
and Steen Pedersen,
{\em  Harmonic analysis and fractal limit-measures induced by
representations of a certain $C\sp *$-algebra},
J. Funct. Anal. {\bf 125} (1994), no.
1, 90--110.

\bibitem[JoPe96]{JoPe96}
Palle E. T. Jorgensen
and Steen Pedersen,
{\em  Harmonic analysis of fractal measures},
Constr. Approx. {\bf 12} (1996), no.
1, 1--30.

\bibitem[JoPe98]{JoPe98}
Palle E. T. Jorgensen
and Steen Pedersen,
{\em Dense analytic
subspaces in fractal $L^{2}$-spaces}, J. Analyse Math. {\bf 75} (1998), 185--228.

\bibitem[JKK05]{JKK05}
  Marius  Junge,
  Peter T. Kim,
and   David W. Kribs,
{\em Universal collective rotation channels and quantum
error correction},
 J. Math. Phys. {\bf 46} (2005), no.
  2, 022102, 18
pp.

\bibitem[Kri05]{Kri05}
  David W.   Kribs,
{\em A quantum computing primer for operator theorists},
 Linear Algebra Appl. {\bf 400}
  (2005), 147--167.

\bibitem[MST99]{MST99}
  Yves  Meyer,
  Fabrice Sellan, and
  Murad S. Taqqu,
{\em Wavelets, generalized white noise and fractional
integration: the synthesis of fractional Brownian motion},
 J. Fourier Anal. Appl. {\bf 5} (1999), no.
  5, 465--494.

\bibitem[Wal23]{Wal23}
J. L. Walsh, {\em A closed set of normal orthogonal functions}, Amer.
J. Math. {\bf 45} (1923), 5--24.

\bibitem[Wic94]{Wic94}
Mladen Victor Wickerhauser, {\em Adapted Wavelet Analysis from
Theory to Software}, A K Peters, Ltd., Wellesley, MA, 1994;
 with a separately available computer disk (IBM-PC or
Macintosh). 

\end{thebibliography}
\end{document}